\documentclass{article}

\usepackage{arxiv}

\usepackage[utf8]{inputenc} 
\usepackage[T1]{fontenc}    
\usepackage{hyperref}       
\usepackage{url}            
\usepackage{booktabs}       
\usepackage{amsfonts}       
\usepackage{nicefrac}       
\usepackage{microtype}      
\usepackage{lipsum}

\usepackage{enumerate}
\usepackage{amsmath}
\usepackage{amssymb}
\usepackage{graphicx}

\newtheorem{theorem}{Theorem}

\newtheorem{proposition}[theorem]{Proposition}

\title{Internet congestion control: from stochastic to dynamical models}

\author{Jos\'{e} M. Amig\'{o}, Angel Gim\'{e}nez, Oscar Mart\'{\i}nez-Bonastre and Jos{\'e} Valero\\
Centro de Investigaci\'{o}n Operativa, \\  Universidad Miguel Hern\'{a}ndez de Elche, \\ Avda. de la Universidad s/n, 03202 Elche, Spain \\
  \texttt{jm.amigo@umh.es, a.gimenez@umh.es, oscar.martinez@umh.es, jvalero@umh.es} 
  \\
}

\begin{document}
\maketitle
\begin{center}
\textit{This paper is dedicated to Manfred Denker on the occasion of his
75th birthday}
\end{center}

\begin{abstract}
Since its inception, control of data congestion on the
Internet has been based on stochastic models. One of the first such models
was Random Early Detection. Later, this model was reformulated as a
dynamical system, with the average queue sizes at a router's buffer being
the states. Recently, the dynamical model has been generalized to improve
global stability. In this paper we review the original stochastic model and
both nonlinear models of Random Early Detection with a two-fold objective:
(i) illustrate how a random model can be \textquotedblleft smoothed
out\textquotedblright\ to a deterministic one through data aggregation, and
(ii) how this translation can shed light into complex processes such as the
Internet data traffic. Furthermore\textbf{, }this paper contains new
materials concerning\textbf{\ }the occurrence of chaos, bifurcation
diagrams, Lyapunov exponents and global stability robustness with respect to
control parameters. The results reviewed and reported here are expected to
help design an active queue management algorithm in real conditions, that
is, when system parameters such as the number of users and the round-trip
time of the data packets change over time. The topic also illustrates the
much-needed synergy of a theoretical approach, practical intuition and
numerical simulations in engineering.
\end{abstract}

\keywords{Internet congestion control \and Adaptive queue management \and
Random early detection \and Discrete-time dynamical systems \and Global stability \and
Robust setting of control parameters}

\textbf{AMS Subject Classification:} 37E05, 39A50, 37G35, 37N35

\section{Introduction}

Internet Protocols (IP) networks include memory buffers to manage the
network traffic by implementing a queue algorithm into buffers. Obviously,
the length of this queue is upper bounded by the size of the buffer, and
they are not effective in their purpose of managing bursts of packets if
they run at full capacity. Therefore, in a scenario of congestion some
actions have to be taken to avoid overflow or a system collapse. An Active
Queue Management (AQM) is an algorithm acting on the queue to control its
length and, thus, achieve an efficient control of the system congestion.
Such management also enables the transmission control protocols (TCPs) to
share links properly.

The Random Early Detection (RED) algorithm \cite{Floyd1993} was the first
formal and complete AQM system implemented in TCP/IP networks \cite{Adams2013}. One of the purposes of designing the RED queue management
algorithm was to detect the beginning of traffic congestion well in advance.
This would allow the system to act in two ways: starting to gradually
discard packets and reducing the transmission rate to that node, both
actions aimed at avoiding scenarios of massive packet discarding, system
collapse or overflow.

The congestion detection indicator used by the RED queue management is an
Exponentially Weighted Moving Average (EWMA) of the queue length. EWMA
actually acts as a low-pass filter that smooths the instantaneous queue
length bursts. The degree of smoothing is determined by a weight $w$. The
algorithm uses a decision mechanism based on a linear probability
distribution to determine when packets are rejected. When a packet reaches
the queue, the algorithm proceeds as follows. If the weighted average of the
queue length is less than a minimum threshold $q_{\min }$, no action is
taken and the packet simply remains in the queue. An average queue length
between $q_{\min }$ and another higher threshold $q_{\max }$ is interpreted
as congestion, in which case an early drop command is executed. More
specifically, the RED algorithm drops an incoming packet with a probability
that depends on the EWMA. This means that adequate $q_{\min }$ and $q_{\max
} $ values are a requirement for good performance. Finally, if the average
queue length is greater than the maximum threshold $q_{\max }$, then it is
understood that the congestion is persistent, and packets are dropped to
avoid a lasting full queue. In this case, the packets are dropped with
probability 1.

RED has been the most studied AQM to date, and has been the basis for the
development of new AQM systems. This is so not only because it was the first
one developed in the Internet community, but also because of several
drawbacks that this algorithm entails, some of which have not yet been
successfully resolved. In particular, it is difficult to adjust the RED
parameters for adequate performance due to their high sensitivity. Moreover,
even if the RED parameters are adjusted properly, they are also very
sensitive to the network conditions. Thus, a set of parameters can work
perfectly for a given load of data traffic but not be suitable if the load
varies slightly. Obviously, this is not a desirable feature since Internet
traffic conditions change rapidly. This is, by far, what has most affected a
widespread implementation of RED. There are many proposals in the literature
to overcome these difficulties via modifications of the RED algorithm, e.g., 
\cite{Arpaci2000,Athuraliya2001,Brandauer2001,Floyd2001,Claypool2004,ChinFu2005}.

Numerical experiments and simulations have usually been the main tool for
setting parameters, also in the papers cited above. As a result, many of the
conclusions published in the literature are based only on numerical results,
thus lacking the necessary mathematical rigor to support them. The first
theoretical study of RED we are aware of was done by Ranjan et al. \cite{Ranjan2004}. This was possible because these authors reformulated the
original, stochastic RED model as a discrete-time dynamical model, thus
allowing them to apply the powerful methods of nonlinear dynamics. In
particular, the authors showed that their dynamical RED model was chaotic.
This first dynamical model was generalized in \cite{Duran2018,Duran2019,Amigo2020} with the objective of enhancing the
controllability of the model. To this end, two new control parameters where
introduced via the probability law for dropping incoming packets. Moreover,
the analysis of the generalized model in \cite{Amigo2020} focused on global
stability rather than local stability like in \cite{Ranjan2004}. Numerical
results (e.g., bifurcation diagrams and parametric domains of robustness)
suggest that, as expected, the generalized model is more stable than the
original one. This is a welcome feature for the prospective design of an
AQM\ in real conditions, that is, when system parameters such as the number
of users and the round-trip time change over time.

In this paper we review the original RED dynamical model \cite{Ranjan2004}
and the main results of the generalized dynamical model \cite{Amigo2020}
with a two-fold objective: (i) illustrate how a random model can be
\textquotedblleft smoothed out\textquotedblright\ to a deterministic one
through data aggregation (sometimes with much ingenuity), and (ii) how this
translation can provide useful insights into complex processes such as data
traffic on the Internet. Furthermore, this paper contains new materials
concerning\textbf{\ }the occurrence of chaos, bifurcation diagrams, Lyapunov
exponents and global stability robustness with respect to control
parameters, as we detail in short. The topic (Internet congestion control)
illustrates also the much-needed synergy of a theoretical approach,
practical intuition and numerical simulations in engineering.

The rest of this paper is organized as follows. In Section 2 we present the
probabilistic RED algorithm and how to derive a nonlinear model by trading
off detailed information (queue size) for coarse-grained information
(average queue length). Our presentation is based on \cite{Ranjan2004}. In
Section 3, we follow \cite{Amigo2020} to generalize the said nonlinear model
with the purpose of increasing the number of control parameters and thus the
control leverage. We have also added in Section 3 the proof that the
generalized RED dynamics can be chaotic (depending on the parameters) in the
sense of Li-Yorke. The main results in \cite{Amigo2020} concerning the
global stability of the unique fixed point of the RED dynamic are summarized
in Section 4, while Section 5 recaps a few more particular results with a
potential to be implemented in real congestion control. The graphics and
numerical simulations presented in Section 6, all of them new, deal
principally with bifurcation diagrams and Lyapunov exponents (Section 6.1),
as well as\ with the robustness of the fixed point with respect to the new
control parameters (Section 6.2). Section 7 contains the concluding remarks
and an outlook.

\section{Random Early Detection}

We follow \cite{Ranjan2004} for the description of RED and the derivation of
its discrete-time dynamical model. Our main concern is the mathematical
content of the model, hence technical details will be kept to a minimum.

Consider the following communication network throughout. $N$ users (i.e., $N$
uniform TCP connections) are connected to a Router 1 which shares an
Internet link (bottleneck link) with Router 2. The \textit{capacity} of this
channel (i.e., the maximum amount of error-free information that can be
transmitted over the channel per unit time) is $C$. Further parameters of
the system are the \textit{packet size} $M$, the \textit{round-trip time}
(propagation delay with no queuing delay) $d$ of the packets, and the 
\textit{buffer size} $B$ of Router 1. So, $B$ is the maximum number of
packets that the buffer of Router 1 can store. The \textit{queue length} or 
\textit{size }$q$ at the buffer at a given time is the number of packets in
the buffer at that time, so $q=0,1,...,B$.

In the RED model, the probability $p$ of dropping an incoming packet at a
router depends on the `average' queue length $q^{\mathtt{ave}}$ (to be
defined below) as follows: 
\begin{equation}
p(q^{\mathtt{ave}})=\left\{ 
\begin{array}{cl}
0 & \text{if }q^{\mathtt{ave}}<q_{\min }, \\ 
1 & \text{if }q^{\mathtt{ave}}>q_{\max }, \\ 
\frac{q^{\mathtt{ave}}-q_{\min }}{q_{\max }-q_{\min }}p_{\max } & \text{otherwise.}\end{array}\right.  \label{p-ave}
\end{equation}Thus, $q_{\min }$ and $q_{\max }$ are the lower and upper thresholds of $q^{\mathtt{ave}}$ for accepting and dropping an incoming packet, respectively,
and $p_{\max }$ is the selected drop probability when $q^{\mathtt{ave}}=q_{\max }$, i.e., the maximum packet drop probability. The average queue
length is updated at the time of a packet arrival according to the averaging\begin{equation}
q_{\mathtt{new}}^{\mathtt{ave}}=(1-w)q_{\mathtt{old}}^{\mathtt{ave}}+wq^{\mathtt{cur}}  \label{q-new}
\end{equation}between the previous average queue length $q_{\mathtt{old}}^{\mathtt{ave}}$
and the current queue length $q^{\mathtt{cur}}$, where $0<w<1$ is the 
\textit{averaging weight}. The initial average queue length is set in an
unspecified way; from then on, the average queue length is defined by the
rule (\ref{q-new}). The higher $w$, the faster the RED mechanism reacts to
the actual buffer occupancy. In practice $w$ is taken rather small,
typically $\lesssim 0.2$ \cite{Wang2014}.

On the way to a discrete-time nonlinear dynamical model, let $q_{n}^{\mathtt{ave}}$ be the average queue length at time $n=1,2,...$ Then, the packet drop
probability at time $n$, $p_{n}$, determines the queue length at time $n+1$, 
$q_{n+1}$. In turn, $q_{n+1}$ is used to compute $q_{n+1}^{\mathtt{ave}}$
according to the averaging rule (\ref{q-new}). Finally, $q_{n+1}^{\mathtt{ave}}$ determines the drop probability at time $n+1$, $p_{n+1}$, via Equation (\ref{p-ave}). In sum, the update of $p_{n}$ has this formal structure:\begin{equation}
\left\{ 
\begin{array}{l}
q_{n+1}=G(p_{n}) \\ 
q_{n+1}^{\mathtt{ave}}=(1-w)q_{n}^{\mathtt{ave}}+wq_{n+1} \\ 
p_{n+1}=H(q_{n+1}^{\mathtt{ave}})\end{array}\right.  \label{structure}
\end{equation}where $H(q_{n+1}^{\mathtt{ave}})=p(q_{n+1}^{\mathtt{ave}})$ is given in (\ref{p-ave}). From (\ref{structure}) it follows\begin{equation}
q_{n+1}^{\mathtt{ave}}=(1-w)q_{n}^{\mathtt{ave}}+wG(H(q_{n}^{\mathtt{ave}})).
\label{master}
\end{equation}So this equation will define a dynamic model for the average queue length
once the function $G$ is determined.

It was shown in \cite{Hespanha2001,Mathis1997,Padhye2000,Ranjan2002} that,
if the packet size $M$, the time delay $d$ and the packet loss probability $p $ are known, the stationary throughput of a typical TCP connection can be
approximated by\begin{equation}
T(p,d)=\frac{MK}{\sqrt{p}d}+o\left( \frac{1}{\sqrt{p}}\right) ,
\label{throughput}
\end{equation}where $K$ is a constant usually set equal to $\sqrt{3/2}$. Therefore, the 
\textit{steady-state packet drop probability} $p_{r}$ such that the
aggregate throughput of connections equals the link capacity, is given by 
\begin{equation}
N\frac{MK}{\sqrt{p_{r}}d}=C.  \label{throughput2}
\end{equation}This is the smallest drop probability such that $p_{n}\geq p_{r}$ entails $q_{n+1}=0$. By (\ref{throughput2}),\begin{equation}
p_{r}=\left( \frac{NMK}{Cd}\right) ^{2}.  \label{p_u}
\end{equation}The corresponding average queue length $q_{r}^{\mathtt{ave}}$ such that $q_{n}^{\mathtt{ave}}>q_{r}^{\mathtt{ave}}$ implies $q_{n+1}=0$, follows from
(\ref{p-ave}):\begin{equation}
q_{r}^{\mathtt{ave}}=\left\{ 
\begin{array}{ll}
\frac{q_{\max }-q_{\min }}{p_{\max }}p_{r}+q_{\min } & \text{if }p_{r}\leq
p_{\max }, \\ 
q_{\max } & \text{otherwise.}\end{array}\right.  \label{q_u}
\end{equation}Note that\begin{equation}
q_{r}^{\mathtt{ave}}\leq q_{\max }\text{.}  \label{q_u2}
\end{equation}

Otherwise, if $p_{n}<p_{r}$ (i.e., $q_{n+1}>0$) and $q_{n+1}\leq B$, the
Internet link capacity between Routers 1 and 2 is fully utilized and the
round-trip time of a packet arriving at time $n+1$ is augmented from $d$ (no
queuing delay) to $d+Mq_{n+1}/C$. From (see (\ref{throughput2})) 
\begin{equation}
N\frac{MK}{\sqrt{p_{n}}(d+\frac{Mq_{n+1}}{C})}=C,  \label{throughput3}
\end{equation}we obtain\begin{equation}
q_{n+1}=\frac{C}{M}\left( \frac{NKM}{\sqrt{p_{n}}C}-d\right) .
\label{throughput4}
\end{equation}Since $q_{n+1}$ is a strictly decreasing function of $p_{n}$, there exists $\min \{p_{n},\,n\geq 1:q_{n+1}=B\}=:p_{l}$, where 
\begin{equation}
p_{l}=\left( \frac{NMK}{BM+Cd}\right) ^{2}  \label{p_l}
\end{equation}by (\ref{throughput4}). The corresponding average queue length is\begin{equation}
q_{l}^{\mathtt{ave}}=\frac{q_{\max }-q_{\min }}{p_{\max }}p_{l}+q_{\min }=\frac{q_{\max }-q_{\min }}{p_{\max }}\left( \frac{NMK}{BM+Cd}\right)
^{2}+q_{\min }.  \label{q_l}
\end{equation}Therefore, $q_{n+1}=B$ for $p_{n}\leq p_{l}$, i.e., $q_{n}^{\mathtt{ave}}\leq q_{l}^{\mathtt{ave}}$. Note that\begin{equation}
q_{l}^{\mathtt{ave}}>q_{\min }\text{.}  \label{q_l2}
\end{equation}

Comparison of (\ref{q_l}) with (\ref{q_u}) shows that\begin{equation}
q_{l}^{\mathtt{ave}}<q_{r}^{\mathtt{ave}}  \label{q_l3}
\end{equation}if $p_{r}<p_{\max }$. For $q_{l}^{\mathtt{ave}}<q_{r}^{\mathtt{ave}}$ to
hold also when $p_{r}\geq p_{\max }$ (i.e., $q_{r}^{\mathtt{ave}}=q_{\max }$), it is necessary that\begin{equation}
\left( \frac{NMK}{BM+Cd}\right) ^{2}<p_{\max }.  \label{constraint consts}
\end{equation}

From Equations (\ref{throughput4}) and the definition of $p_{r}$, Equation (\ref{p_u}), and $p_{l}$, Equation (\ref{q_l}), it follows that the function $G$ in (\ref{structure}) is defined as\begin{equation}
q_{n+1}=G(p_{n})=\begin{cases}
B & \text{if }p_{n}\leq p_{l}, \\ 
\frac{NK}{\sqrt{p_{n}}}-\frac{Cd}{M} & \text{if }p_{l}<p_{n}<p_{r}, \\ 
0 & \text{if }p_{n}\geq p_{r}.\end{cases}
\label{G2}
\end{equation}Finally, plugging (\ref{G2}) into (\ref{master}) we conclude that the time
evolution of the average queue length is given by the nonlinear dynamic\begin{equation}
q_{n+1}^{\mathtt{ave}}=\begin{cases}
(1-w)q_{n}^{\mathtt{ave}}+wB & \text{if }q_{n}^{\mathtt{ave}}\leq q_{l}^{\mathtt{ave}}, \\ 
(1-w)q_{n}^{\mathtt{ave}}+w\left( \frac{NK}{\sqrt{p_{n}}}-\frac{Cd}{M}\right)
& \text{if }q_{l}^{\mathtt{ave}}<q_{n}^{\mathtt{ave}}<q_{r}^{\mathtt{ave}},
\\ 
(1-w)q_{n}^{\mathtt{ave}} & \text{if }q_{n}^{\mathtt{ave}}\geq q_{r}^{\mathtt{ave}},\end{cases}
\label{Ranjan}
\end{equation}where, see (\ref{p-ave}),\begin{equation}
p_{n}=p(q_{n}^{\mathtt{ave}})=\frac{q_{n}^{\mathtt{ave}}-q_{\min }}{q_{\max
}-q_{\min }}p_{\max }  \label{p_n}
\end{equation}since $q_{l}^{\mathtt{ave}}>q_{\min }$ by (\ref{q_l2}) and $q_{r}^{\mathtt{ave}}\leq q_{\max }$ by (\ref{q_u2}).

The mapping $q_{n}^{\mathtt{ave}}\mapsto q_{n+1}^{\mathtt{ave}}$ of the
interval $[0,B]$, Equation (\ref{Ranjan}), is affine on $[0,q_{l}^{\mathtt{ave}}]$, $\cup $-convex on $[q_{l}^{\mathtt{ave}},q_{r}^{\mathtt{ave}}]$ and
linear on $[q_{r}^{\mathtt{ave}},B]$, where $[q_{l}^{\mathtt{ave}},q_{r}^{\mathtt{ave}}]$ can be shown to be a trapping interval of the ensuing
dynamic. It is, therefore, not surprising that this dynamic may be chaotic,
depending on the parameter settings \cite{Ranjan2004}. This is bad news for
the sake of a stable control of the Internet data traffic, which is our main
concern. The RED dynamical model (\ref{Ranjan}) has 5 \textit{system
parameters} ($N,C,d,M,B$) and 4 \textit{user parameters} ($p_{\max },q_{\min
},q_{\max },w$). The former are constrained by condition (\ref{constraint
consts}). The latter are also called \textit{control parameters} because
they can be tuned at will to stabilize the dynamic if necessary. Among the
system parameters, the most critical ones are $N$ and $d$ because they may
change in real time. When it comes to the implementation of adaptive
control, the averaging weight $w$ is the most practical choice. In the next
section we will add to $w$ two more control parameters to improve the
stability properties of the RED dynamical model (\ref{Ranjan}).

\section{A generalized RED dynamical model}

For a more compact notation, we introduce normalized state variables and
thresholds,\begin{equation}
x_{n}=\frac{q_{n}^{\mathtt{ave}}}{B},\;x_{\min }=\frac{q_{\min }}{B},\;x_{\max }=\frac{q_{\max }}{B},  \label{norm variables}
\end{equation}as well as the dimensionless constants\begin{equation}
A_{1}=\frac{NK}{\sqrt{p_{\max }}B},\;\;A_{2}=\frac{Cd}{MB}.  \label{A12}
\end{equation}Actually, out of the parameters in (\ref{A12}) only $C$, $d$ and $M$ are
dimensional (given, say, in kilobytes per second, seconds and kilobytes,
respectively); the other parameters are pure numbers. The constant $K$ will
be set equal to $\sqrt{3/2}$ for definiteness, see (\ref{throughput}).
According to Equation (\ref{constraint consts}), $A_{1}$ and $A_{1}$ are
related as follows.

\begin{proposition}
\label{Proposition1}The constants $A_{1}$ and $A_{2}$ of the RED model (\ref{Ranjan}) are subject to the constraint\begin{equation}
A_{1}<A_{2}+1.  \label{Prop1}
\end{equation}
\end{proposition}

Inequality (\ref{Prop1}) is assumed throughout this paper. Since $p_{\max }$
appears in the denominator of $A_{1}$, the larger $p_{\max }$, the better
for the constraint (\ref{Prop1}) to be fulfilled.

On the grounds explained in Sections 1 and 2, the RED dynamical model (\ref{Ranjan}) was generalized in \cite{Amigo2020} by replacing the probability
law (\ref{p_n}) by\begin{equation}
p(x_{n})=I_{\alpha ,\beta }(z(x_{n}))\cdot p_{\max },  \label{p_n2}
\end{equation}where $I_{\alpha ,\beta }(z)$, $0\leq z\leq 1$, is the beta distribution
function (or normalized incomplete beta function), 
\begin{equation}
I_{\alpha ,\beta }(z)=\frac{\mathfrak{B}(z;\alpha ,\beta )}{\mathfrak{B}(1;\alpha ,\beta )},\;\;\;\;\mathfrak{B}(z;\alpha ,\beta
)=\int_{0}^{z}t^{\alpha -1}(1-t)^{\beta -1}dt,  \label{I(x)}
\end{equation}with $\alpha ,\beta >0$, and\begin{equation}
z(x_{n})=\frac{x_{n}-x_{\min }}{x_{\max }-x_{\min }},\;\;x_{\min }\leq
x_{n}\leq x_{\max }.  \label{z(q)}
\end{equation}Since $I_{1,1}(z)=z$, we recover the conventional RED model (\ref{Ranjan})
for $\alpha =\beta =1$. The beta distribution is related to the chi-square
distribution \cite{Abramo1972}. By definition, $I_{\alpha ,\beta }(z)$ is
strictly increasing, hence invertible. Its inverse, $I_{\alpha ,\beta
}^{-1}(z)$, is also strictly increasing.

\begin{figure}[tbph]
\centering
\includegraphics[width=0.45\textwidth]{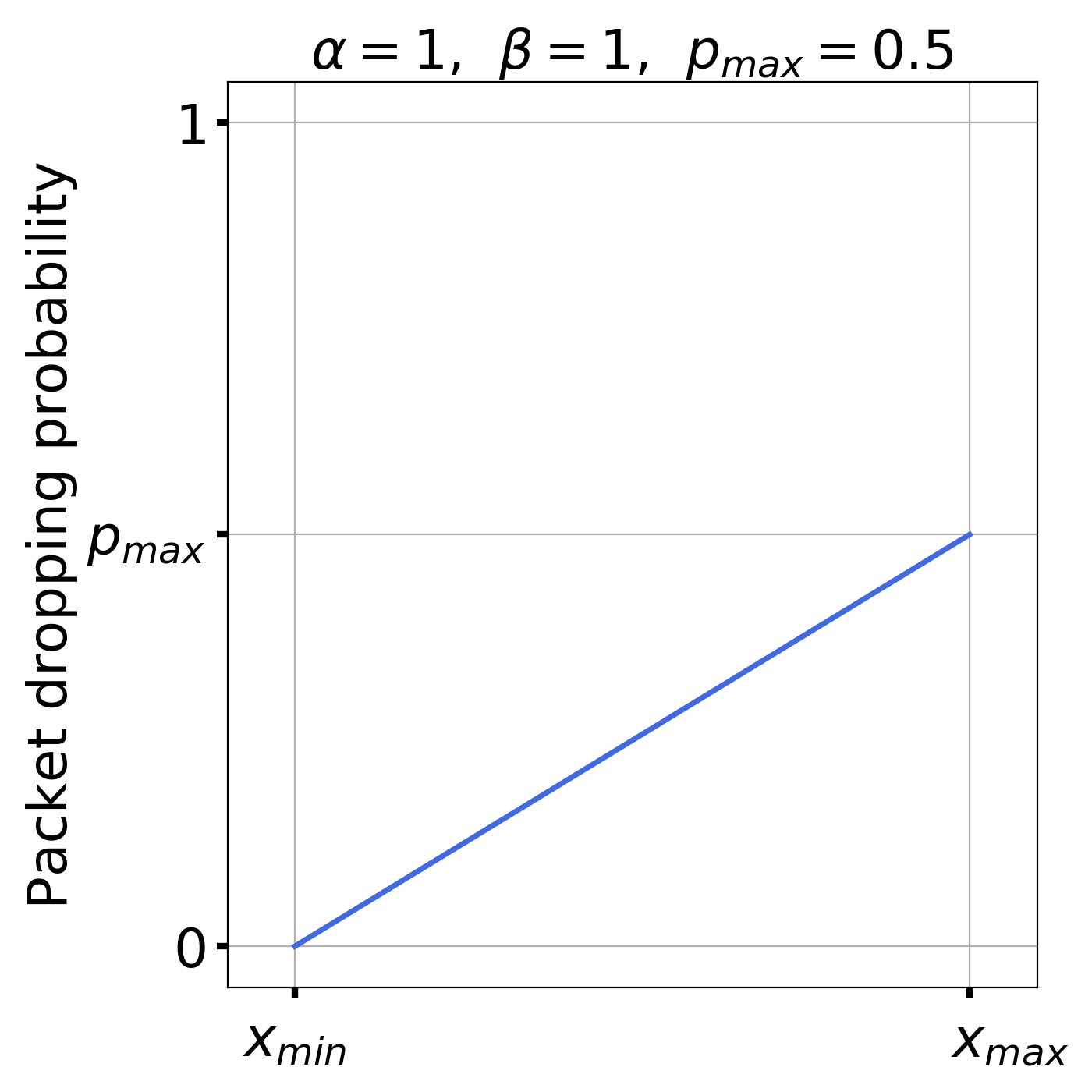} \includegraphics[width=0.45\textwidth]{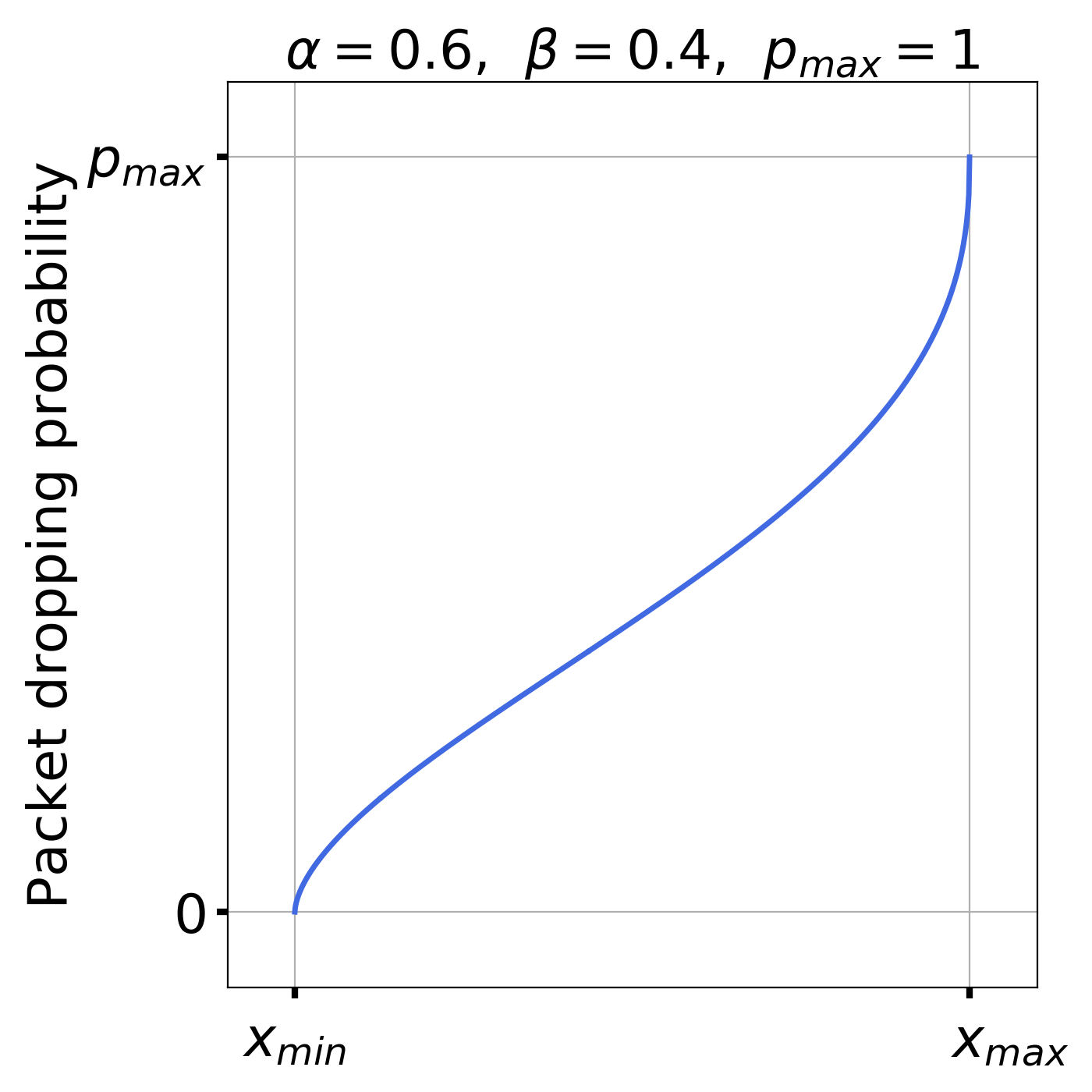}
\caption{Packet dropping probability, Equations (\protect\ref{p_n2})-(\protect\ref{z(q)}), for the settings of the control parameters $\protect\alpha ,\protect\beta $ and $p_{\max }$ displayed at the top of the panels. }
\label{fig-1}
\end{figure}

The packet dropping probability $p(x_{n})$ is plotted in Figure \ref{fig-1} for
different values of the control parameters $\alpha $, $\beta $ and $p_{\max
} $. The left panel ($\alpha =\beta =1$) corresponds to (\ref{p_n}) with $p_{\max }=0.5$.

Thus, we consider hereafter a dynamical system,\begin{equation}
x_{n+1}=f(x_{n}),  \label{dynamics}
\end{equation}where the mapping $f:[0,1]\rightarrow \lbrack 0,1]$ is defined as\begin{equation}
f(x)=\left\{ 
\begin{array}{ll}
(1-w)x+w & \text{if }0\leq x\leq \theta _{l}, \\ 
(1-w)x+w\left( \frac{A_{1}}{\sqrt{I_{\alpha ,\beta }(z(x))}}-A_{2}\right) & 
\text{if }\theta _{l}<x<\theta _{r}, \\ 
(1-w)x & \text{if }\theta _{r}\leq x\leq 1,\end{array}\right.  \label{GDM}
\end{equation}and the thresholds $\theta _{l}$, $\theta _{r}$ are given by

\begin{equation}
\theta _{l}=(x_{\max }-x_{\min })I_{\alpha ,\beta }^{-1}\left( z_{1}\right)
+x_{\min },\text{\ \ }\;z_{1}=\left( \frac{A_{1}}{A_{2}+1}\right) ^{2},
\label{theta_lB}
\end{equation}($z_{1}<1$ by (\ref{constraint consts})) and 
\begin{equation}
\theta _{r}=\left\{ 
\begin{array}{ll}
(x_{\max }-x_{\min })I_{\alpha ,\beta }^{-1}(z_{2})+x_{\min } & \text{if }z_{2}=\left( \frac{A_{1}}{A_{2}}\right) ^{2}\leq 1, \\ 
x_{\max } & \text{otherwise.}\end{array}\right.  \label{theta_rB}
\end{equation}Hence, $\theta _{l}>x_{\min }$, 
\begin{equation}
\left\{ 
\begin{array}{ll}
\theta _{r}<x_{\max } & \text{if }0<A_{1}<A_{2}, \\ 
\theta _{r}=x_{\max } & \text{if }A_{2}\leq A_{1}<A_{2}+1,\end{array}\right.  \label{theta_r <=}
\end{equation}and $\theta _{l}<\theta _{r}$ by (\ref{Prop1}).

The thresholds $\theta _{l}$ and $\theta _{r}$ are set so that $f$ is
continuous on $[0,1]$, except when $A_{1}>A_{2}$, in which case $f$ is lower
semicontinuous at $\theta _{r}$. Indeed, if $A_{1}>A_{2}$, 
\begin{eqnarray}
f(\theta _{r}-) &=&f(x_{\max }-)=\lim_{x\rightarrow x_{\max
}-}f(x)=(1-w)x_{\max }+w(A_{1}-A_{2})  \notag \\
&=&f(x_{\max })+w(A_{1}-A_{2})>f(x_{\max })=f(\theta _{r}),
\label{f(theta_r-)}
\end{eqnarray}where we used $I_{\alpha ,\beta }(z(x_{\max }))=I_{\alpha ,\beta }(1)=1$, on
the second line.

We will numerically show in Section 6 that the additional control parameters 
$\alpha $ and $\beta $ of the generalized RED nonlinear model (\ref{GDM}) do
improve the stability of the dynamic. This is important to enable
controllability also in a real environment, where the system parameters may
change over time. It is worth noting that $p_{\max }$ and all 5 system
parameters go into the dynamic through the constants $A_{1}$ and $A_{2}$,
which amounts to a high degeneracy of the dynamic with respect to the system
parameters. Note also that $A_{1}$, $A_{2}$ and the control parameters
except $w$ go into the thresholds $\theta _{l}$ and $\theta _{r}$ as well.

We show next that the dynamical system defined by (\ref{GDM}) can be
chaotic, depending of the parameter setting. To be more specific, we are
going to prove analytically that the dynamic $x_{n+1}=f(x_{n})$ may exhibit
chaos in the sense of Li-Yorke \cite{Li1975}. Numerical evidence in the form
of bifurcation diagrams and Lyapunov exponents will be presented in Section
6.

We recall the main result of \cite{Li1975} first.

\begin{theorem}
\label{Thm LiYorke}Let $J$ be a closed interval and $F:J\rightarrow J$ a
continuous mapping. Suppose that there is $x_{0}\in J$ such that\begin{equation}
F^{3}(x_{0})\leq x_{0}<F(x_{0})<F^{2}(x_{0})  \label{LY1}
\end{equation}or\begin{equation}
F^{3}(x_{0})\geq x_{0}>F(x_{0})>F^{2}(x_{0}).  \label{LY2}
\end{equation}Then

\begin{description}
\item[(1)] for every $k\in N$ there is a periodic point in $J$ with period $k $;

\item[(2)] there is an uncountable set $S\subset J$ containing no periodic
points, which satisfies the following conditions:

\begin{description}
\item[(i)] For every $x,y\in S$, $x\neq y$,\begin{equation*}
\underset{n\rightarrow \infty }{\lim \sup }\left\vert
F^{n}(x)-F^{n}(y)\right\vert >0
\end{equation*}and\begin{equation*}
\underset{n\rightarrow \infty }{\lim \inf }\left\vert
F^{n}(x)-F^{n}(y)\right\vert =0.
\end{equation*}

\item[(ii)] For every $x\in S$ and periodic point $p\in J,$\begin{equation*}
\underset{n\rightarrow \infty }{\lim \sup }\left\vert
F^{n}(x)-F^{n}(p)\right\vert =0.
\end{equation*}
\end{description}
\end{description}
\end{theorem}

In our case, $J=[0,1]$ and $F=f$, Equation (\ref{GDM}). For $f$ to be
continuous, we assume $A_{1}\leq A_{2}$ for the time being. Following \cite{Ranjan2004}, we choose $x_{0}=\theta _{r}/(1-w)$ with $w<1-\theta _{r}$ (so
that $\theta _{r}<x_{0}<1$). Therefore,\begin{equation*}
f(x_{0})=f\left( \frac{\theta _{r}}{1-w}\right) =\theta _{r},
\end{equation*}and\begin{equation*}
f^{2}(x_{0})=f\left( \theta _{r}\right) =(1-w)\theta _{r},
\end{equation*}hence\begin{equation*}
x_{0}>f(x_{0})>f^{2}(x_{0}).
\end{equation*}By (\ref{LY2}) we only need that $f^{3}(x_{0})\geq x_{0}$ to conclude the
presence of chaos in the form stated in the points (1) and (2) of Theorem \ref{Thm LiYorke}. There are two cases: $f^{2}(x_{0})>\theta _{l}$, i.e., $w<(\theta _{r}-\theta _{l})/\theta _{r}$ (Case I) and $f^{2}(x_{0})\leq
\theta _{l}$, i.e., $w\geq (\theta _{r}-\theta _{l})/\theta _{r}$ (Case II).

In Case I, $\theta _{l}<f^{2}(x_{0})<\theta _{r}$, hence $f^{3}(x_{0})=f(f^{2}(x_{0}))\geq x_{0}$ if and only if\begin{equation}
f\left( (1-w)\theta _{r}\right) =(1-w)^{2}\theta _{r}+w\left( \frac{A_{1}}{\sqrt{I_{\alpha ,\beta }(z((1-w)\theta _{r})}}-A_{2}\right) \geq \frac{\theta _{r}}{1-w},  \label{chaos1}
\end{equation}where $z((1-w)\theta _{r})=((1-w)\theta _{r}-x_{\min })/(x_{\max }-x_{\min
}) $. To verify that inequality (\ref{chaos1}) can be satisfied, use the
following lower bound of the lhs of (\ref{chaos1}):\begin{eqnarray}
f\left( (1-w)\theta _{r}\right) &\geq &(1-w)^{2}\theta _{r}+\frac{wA_{1}}{\sqrt{I_{\alpha ,\beta }(z((1-w)\theta _{r})}}  \notag \\
&\geq &(1-w)^{2}\theta _{r}+\frac{wA_{1}}{\sqrt{I_{\alpha ,\beta }(z(\theta
_{r})}},  \label{lhs}
\end{eqnarray}since $I_{\alpha ,\beta }(z)^{-1/2}$ is a monotone decreasing function. Set
now $w=\epsilon \ll 1$. Then, a sufficient condition for (\ref{chaos1}) to
hold is that the lower bound (\ref{lhs}) is greater than $\frac{\theta _{r}}{1-\epsilon }$, that is, 
\begin{equation*}
(1-2\epsilon )\theta _{r}+\frac{\epsilon A_{1}}{\sqrt{I_{\alpha ,\beta
}(z(\theta _{r})}}\geq (1+\epsilon )\theta _{r}
\end{equation*}to first order in $\epsilon $. It follows that\begin{equation}
\theta _{r}\sqrt{I_{\alpha ,\beta }(z(\theta _{r})}\leq \frac{A_{1}}{3}
\label{chaos1B}
\end{equation}suffices for $f^{3}(x_{0})\geq x_{0}$ to hold, provided $w$ is sufficiently
small.

In Case II, $f^{2}(x_{0})\leq \theta _{l}$, hence $f^{3}(x_{0})=f(f^{2}(x_{0}))\geq x_{0}$ if and only if\begin{equation*}
f\left( (1-w)\theta _{r}\right) =(1-w)^{2}\theta _{r}+w\geq \frac{\theta _{r}}{1-w},
\end{equation*}hence\begin{equation}
\theta _{r}\leq \frac{1-w}{3-w+w^{2}}<\frac{1}{3}.  \label{chaos2}
\end{equation}

In sum, compliance with (\ref{chaos1}) if $w<(\theta _{r}-\theta
_{l})/\theta _{r}$, or with (\ref{chaos2}) otherwise, implies chaos in the
sense of Theorem \ref{Thm LiYorke}. Unlike condition (\ref{chaos1}), which
involves all the parameters of the model, condition (\ref{chaos2}) only
involves the control parameter $w$. Actually, we expect the condition (\ref{chaos1B}) to be the general criterion because, as mentioned in Section 2,
typically $w\lesssim 0.2$. Needless to say, the possibility of a chaotic
behavior poses a challenge to the control design.

\section{Global stability of the generalized RED model: main results}

In this paper we are interested in a stable and robust congestion control of
the Internet data traffic based on the generalized RED model. This means
that the system remains on, or evolves within the basin of attraction of, a
fixed point for suitably chosen control parameters, also when the system
parameters (grouped in the constants $A_{1}$ and $A_{2}$) and the control
parameters (mainly $w$, $\alpha $ and $\beta $) change slightly. This way
one can expect to keep the system in a neighborhood of the fixed point when
both system parameters and (consequently) control parameters change in real
time. Recall that $A_{1}<A_{2}+1$ (Proposition \ref{Proposition1}).

To this end, this section summarizes the most important results on global
stability of the dynamical system (\ref{GDM}). Proofs and additional detail
can be found in \cite{Amigo2020}. The study of robustness with respect to
the control parameters $\alpha $, $\beta $ and $w$ is deferred to Section
6.2.

\begin{theorem}
\label{Thm1}The mapping $f$ has a unique fixed point $x^{\ast }\in (\theta
_{l},\theta _{r})$ if and only if $A_{1}<A_{2}+x_{\max }$; otherwise, it has
none. Furthermore, $x^{\ast }$ does not depend on $w$.
\end{theorem}

To compute the fixed point $x^{\ast }$ one has in general to solve
numerically the equation $f(x)=x$ (with $\theta _{l}<x<\theta _{r}$ and $A_{1}<A_{2}+x_{\max }$), which leads to\begin{equation}
\frac{A_{1}}{\sqrt{I_{\alpha ,\beta }(z(x))}}=x+A_{2},  \label{fixed point}
\end{equation}where the lhs is a strictly decreasing function of $x$ and the rhs is a
strictly increasing function. It follows that $x^{\ast }$ is an increasing
function of $A_{1}$ (in particular, of $N$) by fixed $A_{2}$, whereas it is
a decreasing function of $A_{2}$ (in particular, of $d$) by fixed $A_{1}$.
In some case (e.g., $\alpha =\beta =1$ and $\alpha =2$, $\beta =1)$,
Equation (\ref{fixed point}) can be solved analytically.

Uniqueness of the fixed point $x^{\ast }$ is good news for the control
stability. But for the sake of controllability we need much more: global
attractiveness of $x^{\ast }$.

Let $\mathcal{B}(S,f)$ denote the basin of attraction of a set $S\subset
\lbrack 0,1]$, that is, $\mathcal{B}(S,f)$ consists of all points of $[0,1]$
that asymptotically end up in $S$.

\begin{theorem}
\label{Thm2}If $(\theta _{l},\theta _{r})$ is invariant (i.e., $f((\theta
_{l},\theta _{r}))\subset (\theta _{l},\theta _{r})$), then $\mathcal{B}((\theta _{l},\theta _{r}),f)=[0,1]$.
\end{theorem}

Theorem \ref{Thm2} shows that the interval $[\theta _{l},\theta _{r}]$ contains the
non-wandering set for $f$ when $(\theta _{l},\theta _{r})$ is invariant,
i.e., it is the dynamical core of the RED dynamic. Actually, the only way to
prevent the orbit of $x_{0}\in \lbrack 0,\theta _{l}]\cup \lbrack \theta
_{r},1]$ from getting trapped within $(\theta _{l},\theta _{r})$ is that it
is a preimage of an hypothetical periodic cycle $f(\theta _{l})=\theta _{r}$
and $f(\theta _{r})=\theta _{l}$. It is easy to show that $\theta
_{l}+\theta _{r}\neq 1$ excludes that possibility.

If $\mathcal{B}(x^{\ast },f)=[0,1]$ we say that the fixed point $x^{\ast }$
is a global attractor of $f$, what amounts to $x^{\ast }$ being \textit{globally stable}. Equivalently, we say also that $f$ is globally stable. Let 
$\left. f\right\vert _{(\theta _{l},\theta _{r})}$ denote as usual the
restriction of the mapping $f$ to the interval $(\theta _{l},\theta _{r})$,\begin{equation}
\left. f\right\vert _{(\theta _{l},\theta _{r})}(x)=(1-w)x+w\left( \frac{A_{1}}{\sqrt{I_{\alpha ,\beta }(z(x))}}-A_{2}\right) .  \label{restrict f}
\end{equation}For further reference, note that $\left. f\right\vert _{(\theta _{l},\theta
_{r})}$ is a smooth function because $z(x)>0$ for all $x\in (\theta
_{l},\theta _{r})$ (since $x\geq \theta _{l}>x_{\min }$); in particular,\begin{eqnarray}
\left. f^{\prime }\right\vert _{(\theta _{l},\theta _{r})}(x) &=&1-w\left( 1+\frac{A_{1}}{2(x_{\max }-x_{\min })}\frac{I_{\alpha ,\beta }^{\prime }(z(x))}{I_{\alpha ,\beta }(z(x))^{3/2}}\right)  \notag \\
&<&1-w,  \label{restrict f'}
\end{eqnarray}$0<w<1$, for all system and control parameters. If $A_{1}>A_{2}$ ($f$
discontinuous at $\theta _{r}$) and $\beta <1$, then $f^{\prime }(\theta
_{r}-)=f^{\prime }(x_{\max }-)=-\infty $, so $\left. f^{\prime }\right\vert
_{(\theta _{l},\theta _{r})}$ is upper bounded but not necessarily lower
bounded.

\begin{theorem}
\label{Thm3}Suppose $A_{1}<A_{2}+x_{\max }$ and let $x^{\ast }\in (\theta
_{l},\theta _{r})$ be the unique fixed point of $f$. If

\begin{description}
\item[(i)] $(\theta _{l},\theta _{r})$ is invariant (so that $\left.
f\right\vert _{(\theta _{l},\theta _{r})}$ defines a dynamical system on $(\theta _{l},\theta _{r})$),

\item[(ii)] $\theta _{l}+\theta _{r}\neq 1$ (so that $\{\theta _{l},\theta
_{r}\}$ is not a periodic orbit), and

\item[(iii)] $\mathcal{B}(x^{\ast },\left. f\right\vert _{(\theta
_{l},\theta _{r})})=(\theta _{l},\theta _{r})$ (i.e., $x^{\ast }$ is a
global attractor of $\left. f\right\vert _{(\theta _{l},\theta _{r})}$),

\item then $x^{\ast }$ is a global attractor of $f$.
\end{description}
\end{theorem}

Therefore, to study the global stability of the RED dynamics it suffices to
concentrate on $\left. f\right\vert _{(\theta _{l},\theta _{r})}$, provided
conditions (i) and (ii) of Theorem \ref{Thm3} are satisfied. Practical
considerations advice to limit the stability analysis to monotone and
unimodal mappings over $(\theta _{l},\theta _{r})$. A bimodal $\left.
f\right\vert _{(\theta _{l},\theta _{r})}$ can be seen in \cite{Amigo2020}.

That a monotone $\left. f\right\vert _{(\theta _{l},\theta _{r})}$ is
strictly increasing or strictly decreasing depends on whether $f(\theta
_{l})<f(\theta _{r}-)$ or $f(\theta _{l})>f(\theta _{r}-)$, respectively.
This translates into the following condition on the averaging weight \cite[Proposition 3]{Amigo2020}:\begin{equation}
f(\theta _{l})\gtrless f(\theta _{r}-)\;\;\Leftrightarrow \;\;w\gtrless 
\frac{\theta _{r}-\theta _{l}}{\theta _{r}-\theta _{l}+1-(A_{1}-A_{2})^{+}}.
\label{><}
\end{equation}As usual, the notation $(A_{1}-A_{2})^{+}$ stands for the positive part of $A_{1}-A_{2}$ ($=\max \{A_{1}-A_{2},0\}$), and similarly for other arguments.

\begin{theorem}
\label{Thm4}If $\left. f\right\vert _{(\theta _{l},\theta _{r})}$ is

\begin{description}
\item[(a)] strictly increasing and $A_{1}<A_{2}+x_{\max }$, or

\item[(b)] strictly decreasing with $\left. f^{\prime }\right\vert _{(\theta
_{l},\theta _{r})}(x)>-1$ and 
\begin{equation}
w\leq \min \left\{ \frac{\theta _{r}-\theta _{l}}{1-\theta _{l}},\frac{\theta _{r}-\theta _{l}}{\theta _{r}-(A_{1}-A_{2})^{+}}\right\}
\label{invar monoton}
\end{equation}(where $\theta _{r}=x_{\max }$ if $(A_{1}-A_{2})^{+}>0$, Equation (\ref{theta_r <=})),
\end{description}

then $\mathcal{B}(q^{\ast },\left. f\right\vert _{(\theta _{l},\theta
_{r})})=(\theta _{l},\theta _{r})$.
\end{theorem}

The restriction $\left. f^{\prime }\right\vert _{(\theta _{l},\theta
_{r})}(x)>-1$ in (b) can be replaced by the absence of 2-cycles, also at the
endpoints $\{\theta _{l},\theta _{r}\}$, should $x^{\ast }$ be a global
attractor. The latter can be done, as in Theorem \ref{Thm3}, with the
proviso $\theta _{l}+\theta _{r}\neq 1$. By Sharkovsky's theorem \cite{Sharko1964} applied to the continuous case ($A_{1}\leq A_{2}$), if there
are no periodic orbits of period 2, then there are no periodic orbits of any
period.

Consider next the unimodal case. Specifically we suppose that $\left.
f\right\vert _{(\theta _{l},\theta _{r})}$ has a local extremum at $x_{c}\in
(\theta _{l},\theta _{r})$.

\begin{theorem}
\label{Thm5}Suppose that (i) $w\leq \frac{\theta _{r}-\theta _{l}}{1-\theta
_{l}}$, (ii) $\left. f\right\vert _{(\theta _{l},\theta _{r})}$ has a local
minimum at $x_{c}\leq x^{\ast }$ (so $f^{\prime }(x^{\ast })\geq 0$), and
(iii) $A_{1}\leq A_{2}+x_{\max }$. Then $(\theta _{l},\theta _{r})$ is
invariant and $\mathcal{B}(x^{\ast },\left. f\right\vert _{(\theta
_{l},\theta _{r})})=(\theta _{l},\theta _{r})$.
\end{theorem}

The case $x_{c}>x^{\ast }$ is more involved.

\begin{theorem}
\label{Thm6}The following holds.

\begin{description}
\item[(a)] Suppose that (i) $x^{\ast }>(A_{1}-A_{2})^{+}$, (ii) $\left.
f\right\vert _{(\theta _{l},\theta _{r})}$ has a minimum at $x_{c}>x^{\ast }$
(so $f^{\prime }(x^{\ast })<0$), and (iii) 
\begin{equation}
w\leq \min \left\{ \frac{\theta _{r}-\theta _{l}}{1-\theta _{l}},\,\frac{x^{\ast }-\theta _{l}}{x^{\ast }-(A_{1}-A_{2})^{+}}\right\} .  \label{wThm3}
\end{equation}If $f^{\prime }(x)>-1$ for $\theta _{l}<x<x_{c}$, then $(\theta _{l},\theta
_{r})$ is invariant and $\mathcal{B}(x^{\ast },\left. f\right\vert _{(\theta
_{l},\theta _{r})})=(\theta _{l},\theta _{r})$.

\item[(b)] Part (a) holds also if (i) is replaced by $x^{\ast }\leq $ $(A_{1}-A_{2})^{+}<x_{\max }$, and (iii) is replaced by\begin{equation}
w\leq \frac{\theta _{r}-\theta _{l}}{1-\theta _{l}}.  \label{wThm3(b)}
\end{equation}
\end{description}
\end{theorem}

Similar results can be obtained when $\left. f\right\vert _{(\theta
_{l},\theta _{r})}$ has a local maximum. However, the assumptions in this
case are more restrictive due to (\ref{restrict f'}). Therefore, we will not
consider this option hereafter. The same happens regarding the Schwarzian
derivative\begin{equation}
Sf(x)=\frac{f^{\prime \prime \prime }(x)}{f^{\prime }(x)}-\frac{3}{2}\left( 
\frac{f^{\prime \prime }(x)}{f^{\prime }(x)}\right) ^{2}.  \label{Schwarz}
\end{equation}The condition $\left. \left. Sf\right. \right\vert _{[\theta _{l},\theta
_{r}]}(x)<0$ for all $x\neq x_{c}$, along with (i) $A_{1}\leq A_{2}$ (so
that $f$ is continuous at $\theta _{r})$, (ii) an invariant $[\theta
_{l},\theta _{r}]$, and (iii) $\left\vert f^{\prime }(q^{\ast })\right\vert
\leq 1$, implies $\mathcal{B}(x^{\ast },\left. f\right\vert _{[\theta
_{l},\theta _{r}]})=[\theta _{l},\theta _{r}]$ \cite[Proposition 1]{ElMorshedy2006}, but its implementation is quite involved and restrictive
in parametric space.

\section{Global stability of the generalized RED model: particular results}

When it comes to putting in place a control mechanism of the dynamics (\ref{GDM}) that can be upgraded to an adaptive mechanism in real time,
simplicity and computational speed is a must. This has been our guideline
when selecting Theorems \ref{Thm4} to \ref{Thm6}, from which one derives the sought global
stability of $x^{\ast }$, i.e., $\mathcal{B}(x^{\ast },f)=[0,1]$, under the
following three assumptions:

(i) $(\theta _{l},\theta _{r})$ is invariant,

(ii) $\{\theta _{l},\theta _{r}\}$ is not a periodic orbit, and

(iii) $\left. f\right\vert _{(\theta _{l},\theta _{r})}$ has at most one
minimum.

\noindent It turns out that a good compromise satisfying (i)-(iii) is to
choose $\left. f\right\vert _{(\theta _{l},\theta _{r})}$ $\cup $-convex. In
particular, the convexity of $f$ in the dynamical core $(\theta _{l},\theta
_{r})$ simplifies Theorems \ref{Thm4}(b) and \ref{Thm6}(a) as follows.

\begin{theorem}
\label{Thm7}If $\left. f\right\vert _{(\theta _{l},\theta _{r})}$ is $\cup $-convex, then the hypothesis $\left. f^{\prime }\right\vert _{(\theta
_{l},\theta _{r})}(x)>-1$ in Theorem \ref{Thm1}(b) may be replaced by $f^{\prime }(\theta _{l}+)\geq -1$.
\end{theorem}

\begin{theorem}
\label{Thm8}If $\left. f\right\vert _{(\theta _{l},\theta _{r})}$ is $\cup $-convex, along with $(A_{1}-A_{2})^{+}<x^{\ast }$ and $x_{c}>x^{\ast }$ (so $f^{\prime }(x^{\ast })<0$), then the assumptions $\left. f^{\prime
}\right\vert _{(\theta _{l},q_{c})}(x)>-1$ and (\emph{\ref{wThm3}}) in
Theorem \ref{Thm6}(a) may be replaced by $f^{\prime }(\theta _{l}+)\geq -1$
and\begin{equation}
w\leq \min \left\{ \frac{\theta _{r}-\theta _{l}}{1-\theta _{l}},\frac{x^{\ast }-\theta _{l}+\frac{1}{m}(\theta _{l}-(A_{1}-A_{2})^{+})^{+}}{x^{\ast }-(A_{1}-A_{2})^{+}}\right\} ,  \label{wThm8}
\end{equation}respectively, where $m=(1-f^{\prime }(x^{\ast })/w>1$.
\end{theorem}

Since $\frac{1}{m}(\theta _{l}-(A_{1}-A_{2})^{+})^{+}\geq 0$, the bound (\ref{wThm8}) is indeed less restrictive than the bound (\ref{wThm3}). When $m\rightarrow \infty $ in (\ref{wThm8}) we recover (\ref{wThm3}), as it
should.

General conditions for $\left. f\right\vert _{(\theta _{l},\theta _{r})}$ to
be $\cup $-convex (with one or no critical point) follow readily from the
expression 
\begin{equation}
\left. f^{\prime \prime }\right\vert _{(\theta _{l},\theta _{r})}(x)=\frac{wA_{1}}{4(x_{\max }-x_{\min })}\frac{I_{\alpha ,\beta }^{\prime }(z(x))}{I_{\alpha ,\beta }(z(x))^{3/2}}\left[ 3J_{\alpha ,\beta }(z(x))-2h_{\alpha
,\beta }(z(x))\right] ,  \label{f''}
\end{equation}where $z(x)=(x-x_{\min })/(x_{\max }-x_{\min })$, and 
\begin{equation}
J_{\alpha ,\beta }(z):=\frac{I_{\alpha ,\beta }^{\prime }(z)}{I_{\alpha
,\beta }(z)}>0,\;\;h_{\alpha ,\beta }(z):=\frac{\alpha -1}{z}-\frac{\beta -1}{1-z},  \label{J(x)}
\end{equation}for all $z(\theta _{l})<z<z(\theta _{r})$.

\begin{proposition}
\label{Proposition2}Suppose that $3J_{\alpha ,\beta }(z(x))-2h_{\alpha
,\beta }(z(x))>0$ for $\theta _{l}<x<\theta _{r}$. It holds:

\begin{description}
\item[(a)] $f_{(\theta _{l},\theta _{r})}$ is $\cup $-convex.

\item[(b)] If $f^{\prime }(\theta _{l}+)\cdot f^{\prime }(\theta _{r}-)>0$,
then $\left. f\right\vert _{(\theta _{l},\theta _{r})}$ has no critical
point; if $f^{\prime }(\theta _{l}+)\cdot f^{\prime }(\theta _{r}-)<0$, then 
$\left. f\right\vert _{(\theta _{l},\theta _{r})}$ has one critical point.
\end{description}
\end{proposition}

There are a number of settings for the control parameters $\alpha $ and $\beta $ that guarantee the sufficient condition $3J_{\alpha ,\beta
}(z(x))-2h_{\alpha ,\beta }(z(x))>0$ for $f_{(\theta _{l},\theta _{r})}$ to
be $\cup $-convex. The perhaps most useful ones are the following.

\begin{proposition}
\label{Proposition3}$\left. f\right\vert _{(\theta _{l},\theta _{r})}$ is $\cup $-convex at $x$ in the following cases:

\begin{description}
\item[(a)] $\alpha \leq \beta $ and $0<z(x)<\frac{\alpha }{\alpha +\beta }$.

\item[(b)] $\alpha >\beta $ and $0<z(x)<\frac{\alpha +2}{\alpha +\beta +4}$.
\end{description}
\end{proposition}

Therefore, the mapping $\left. f\right\vert _{(\theta _{l},\theta _{r})}$ is 
$\cup $-convex whenever 
\begin{equation*}
z(\theta _{r})=\frac{\theta _{r}-x_{\min }}{x_{\max }-x_{\min }}\leq \left\{ 
\begin{array}{cl}
\frac{\alpha }{\alpha +\beta } & \text{if }\alpha \leq \beta , \\ 
\frac{\alpha +2}{\alpha +\beta +4} & \text{if }\alpha >\beta .\end{array}\right.
\end{equation*}

\section{Numerical simulations}

The system parameters used as reference values in this section are\begin{equation}
N=1850,\;C=321,000\text{ kBps},\;d=0.012\text{ s},\;M=1\text{ kB},\;B=2000\text{ packets},  \label{UMH}
\end{equation}along with $K=\sqrt{3/2}\simeq 1.225$. Data (\ref{UMH}) correspond to the
Miguel Hern\'{a}ndez University network and are the same as in \cite{Amigo2020}. The parameters $A_{1}$ and $A_{2}$ result then in 
\begin{equation*}
A_{1}\simeq \frac{1.13}{\sqrt{p_{\max }}}\text{\ and\ \ }A_{2}\simeq 1.93,
\end{equation*}see (\ref{A12}), so the constraint $A_{1}<A_{2}+1$, Proposition \ref{Proposition1}, holds for $p_{\max }\gtrsim 0.15$.

As for the control parameters: $\alpha $, $\beta $ and $p_{\max }$ will be
specified in each figure, while the reference values for the remaining
control parameters are 
\begin{equation}
w=0.15,\;\;x_{\min }=0.2,\;\;x_{\max }=0.6  \label{UMH2}
\end{equation}when kept fixed.

Next we present bifurcation diagrams of the RED dynamics and the
corresponding Lyapunov exponents (Section 6.1), as well as a brief
discussion of the robustness of the fixed point $x^{\ast }$ with respect to
the control parameters $\alpha $, $\beta $ (Section 6.2). The primary
objective of Section 6.1 is to illustrate the effect of the control
parameters $\alpha ,\beta $ on two particular bifurcation diagrams. As for
Section 6.2, robustness of the RED dynamics with respect to the control
parameters in general, and with respect to $\alpha $ and $\beta $ in
particular, is essential to enable controllability under real conditions.

\subsection{Bifurcation diagrams and Lyapunov exponents}

\begin{figure}[tbph]
\centering
\includegraphics[width=0.45\textwidth]{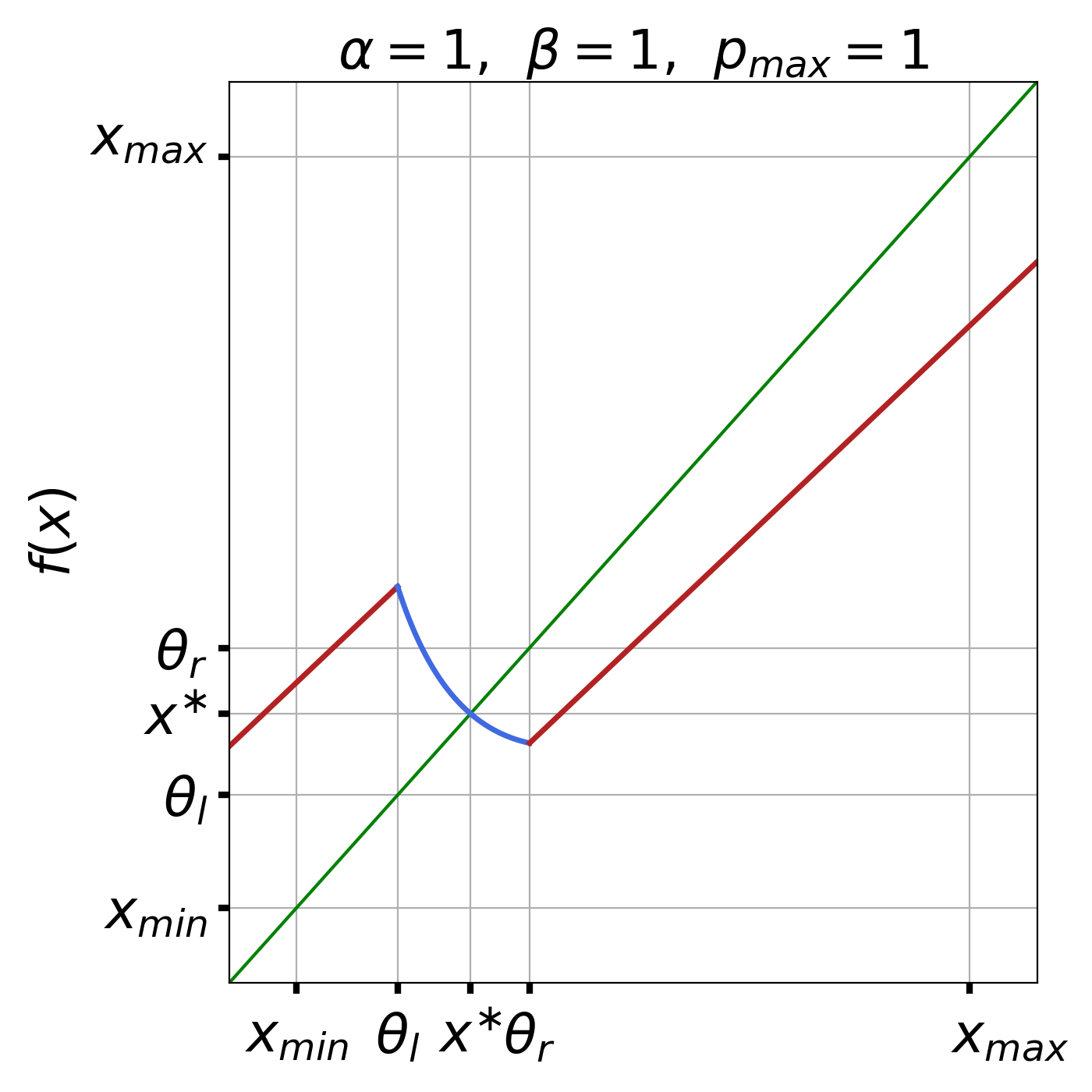} \includegraphics[width=0.45\textwidth]{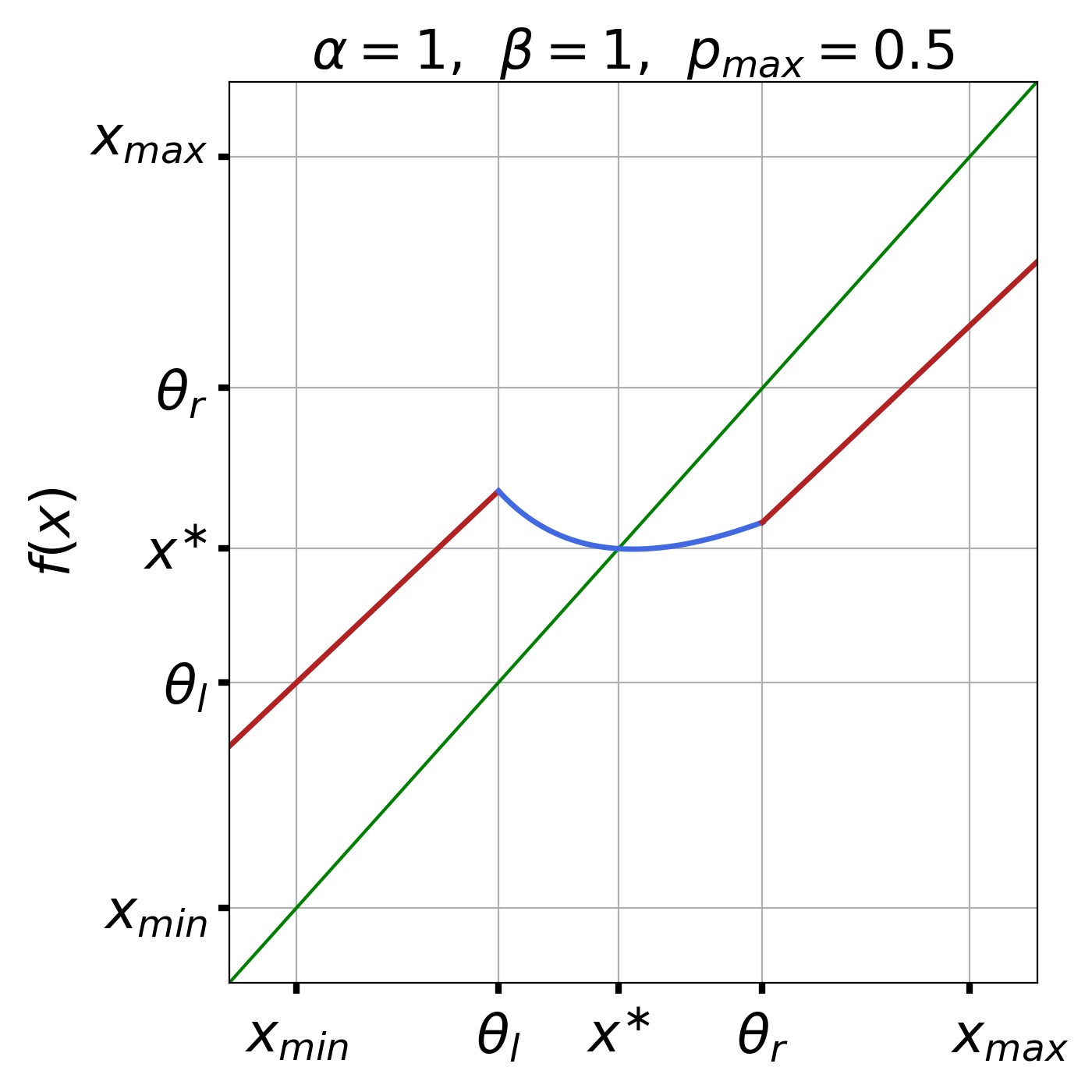} \includegraphics[width=0.45\textwidth]{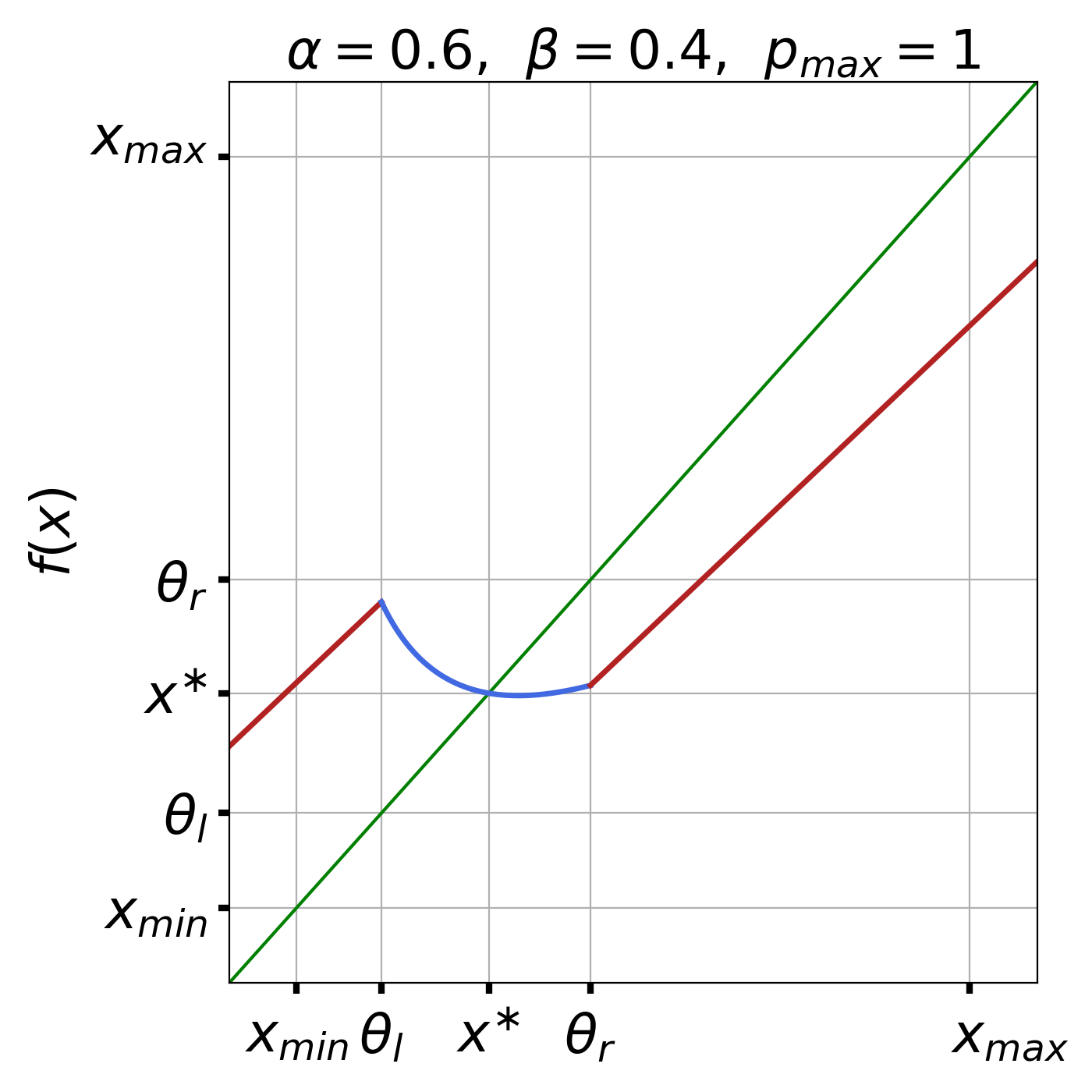} \includegraphics[width=0.45\textwidth]{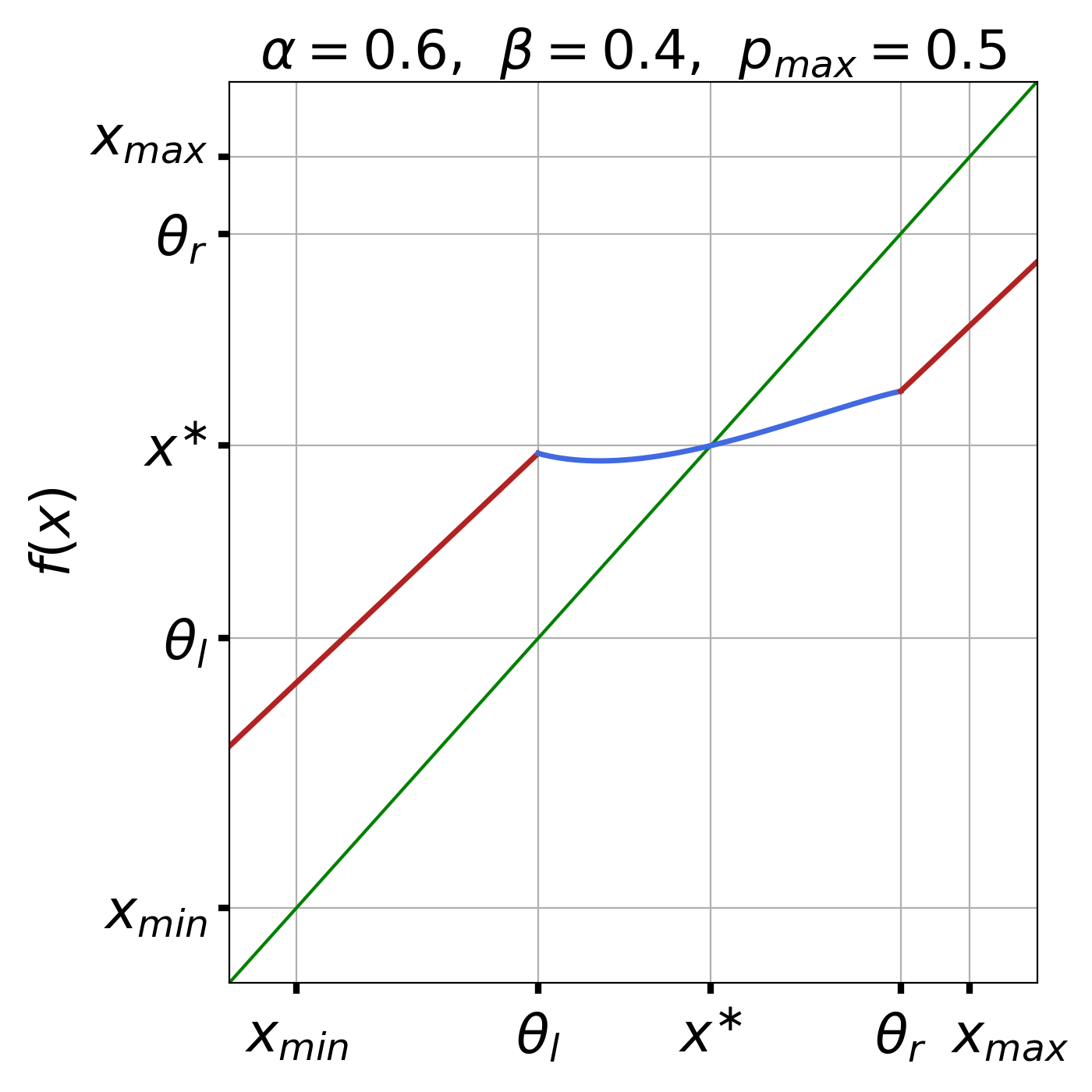}
\caption{Graphs of the mapping $f$, Equation (\protect\ref{GDM}), for the
values of $\protect\alpha $, $\protect\beta $ and $p_{\max }$ displayed at
the top of the panels. Remaining parameters are given in (\protect\ref{UMH})
and (\protect\ref{UMH2}). The restriction $\left. f\right\vert _{(\protect\theta _{l},\protect\theta _{r})}$ is monotonically decreasing in the upper
left panel, while it is $\cup $-convex in the other panels. }
\label{fig-2}
\end{figure}

The choice of appropriate parameters to ensure the stability of the dynamic
(meaning that the fixed point $x^{\ast }$ is a global attractor) is perhaps
the biggest drawback of the RED algorithm. Therefore, it is of great
interest to analyze the bifurcation diagram and the Lyapunov exponent for
different parameters. Some of them were studied in \cite{Amigo2020}, but the
analysis of the two very important parameters $x_{\min }$ and $x_{\max }$
were omitted there for brevity. The mappings we use here are plotted in
Figure \ref{fig-2}.

\begin{figure}[tbph]
\centering
\includegraphics[width=0.45\textwidth]{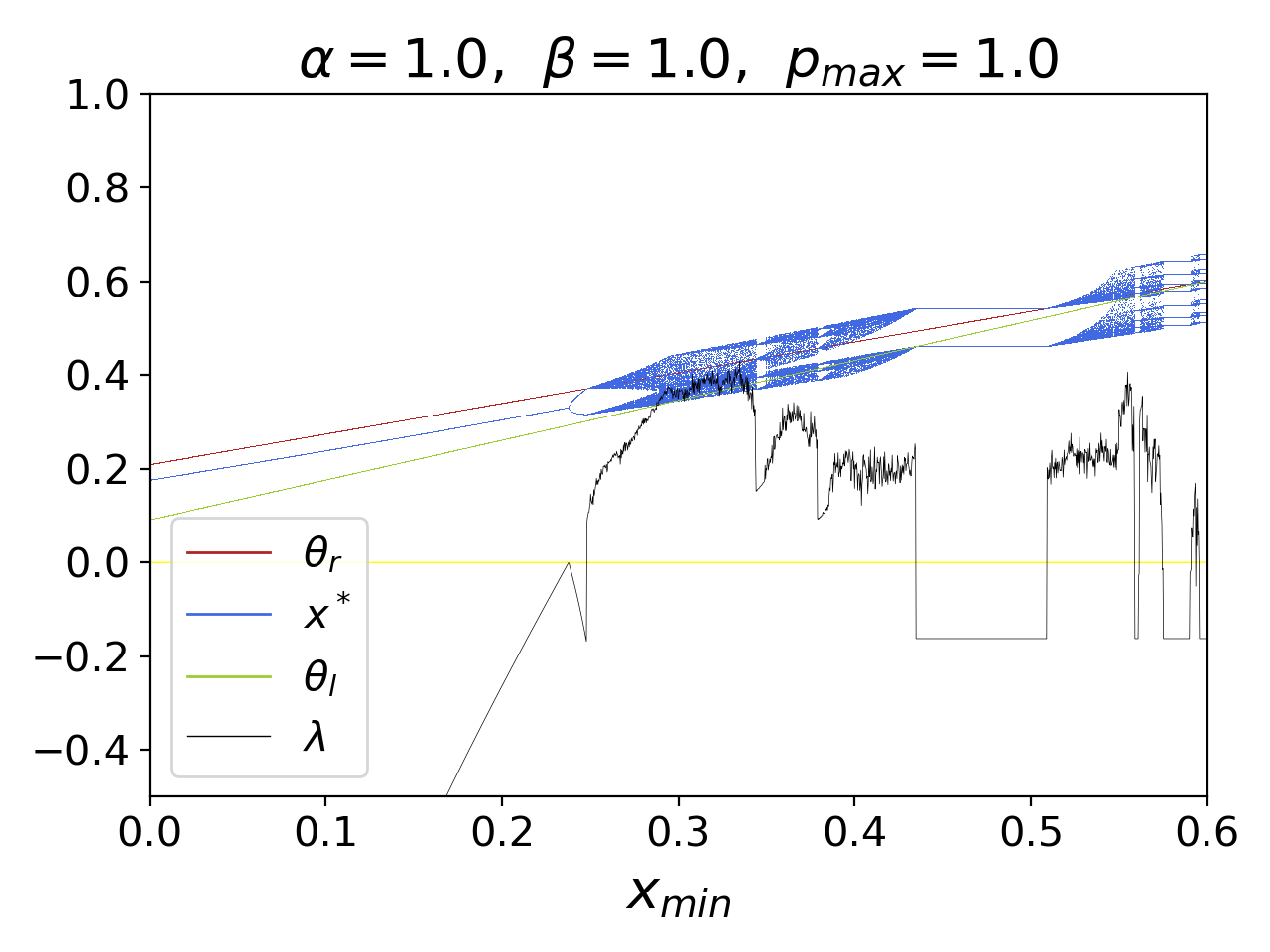} \includegraphics[width=0.45\textwidth]{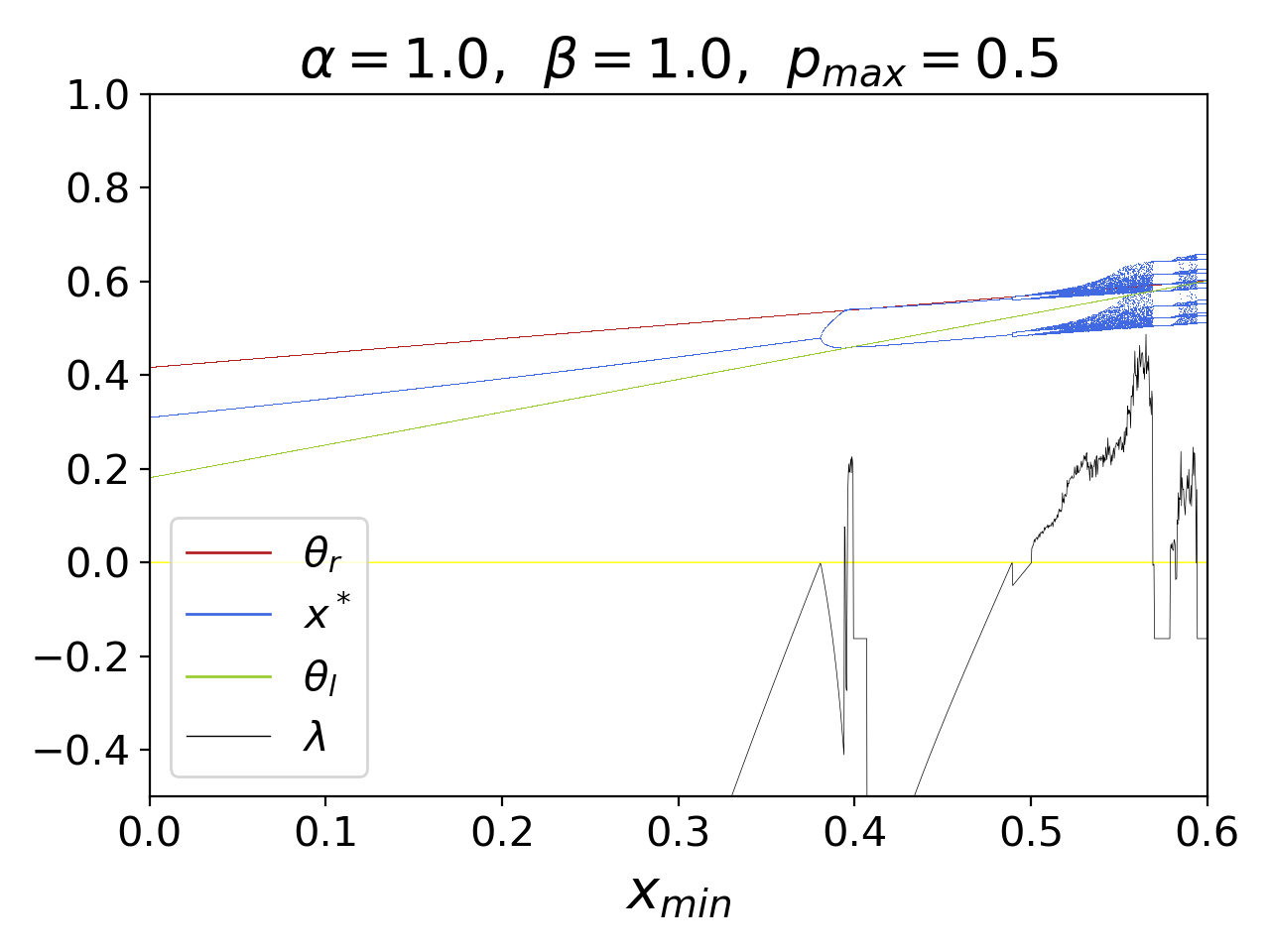} \includegraphics[width=0.45\textwidth]{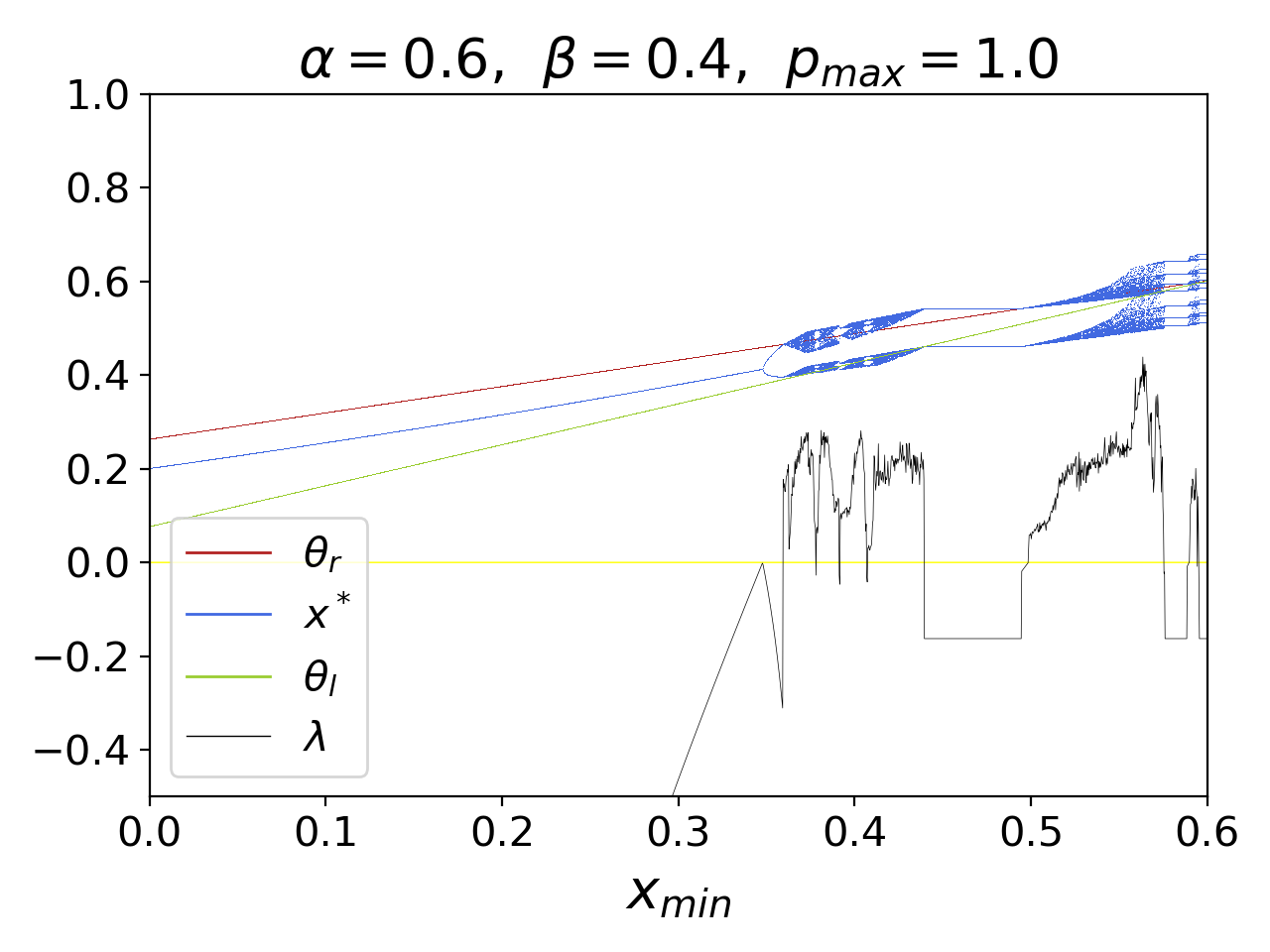} \includegraphics[width=0.45\textwidth]{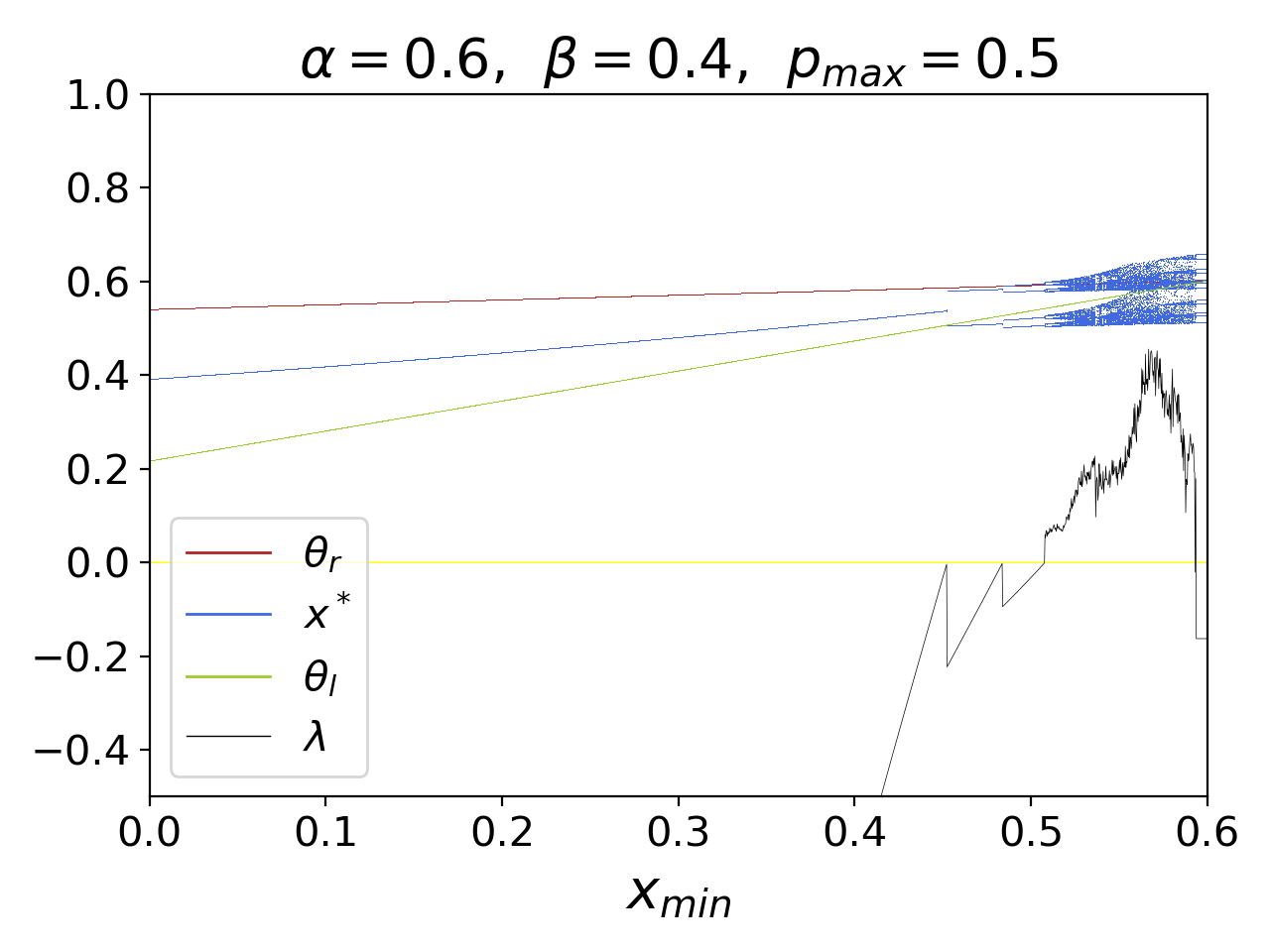}
\caption{Bifurcation diagram and Lyapunov exponents of the normalized
average queue length $x$ with respect to the parameter $x_{\min }$ for the
values of $\protect\alpha $, $\protect\beta $ and $p_{\max }$ displayed at
the top of the panels. The parameter $x_{\min }$ ranges in the interval $[0.2,0.6)$. System parameters are given in (\protect\ref{UMH}). Other
control parameters: $x_{\max }=0.6$ and $w=0.15$.}
\label{fig-bif-xmin}
\end{figure}

\begin{figure}[tbph]
\centering
\includegraphics[width=0.45\textwidth]{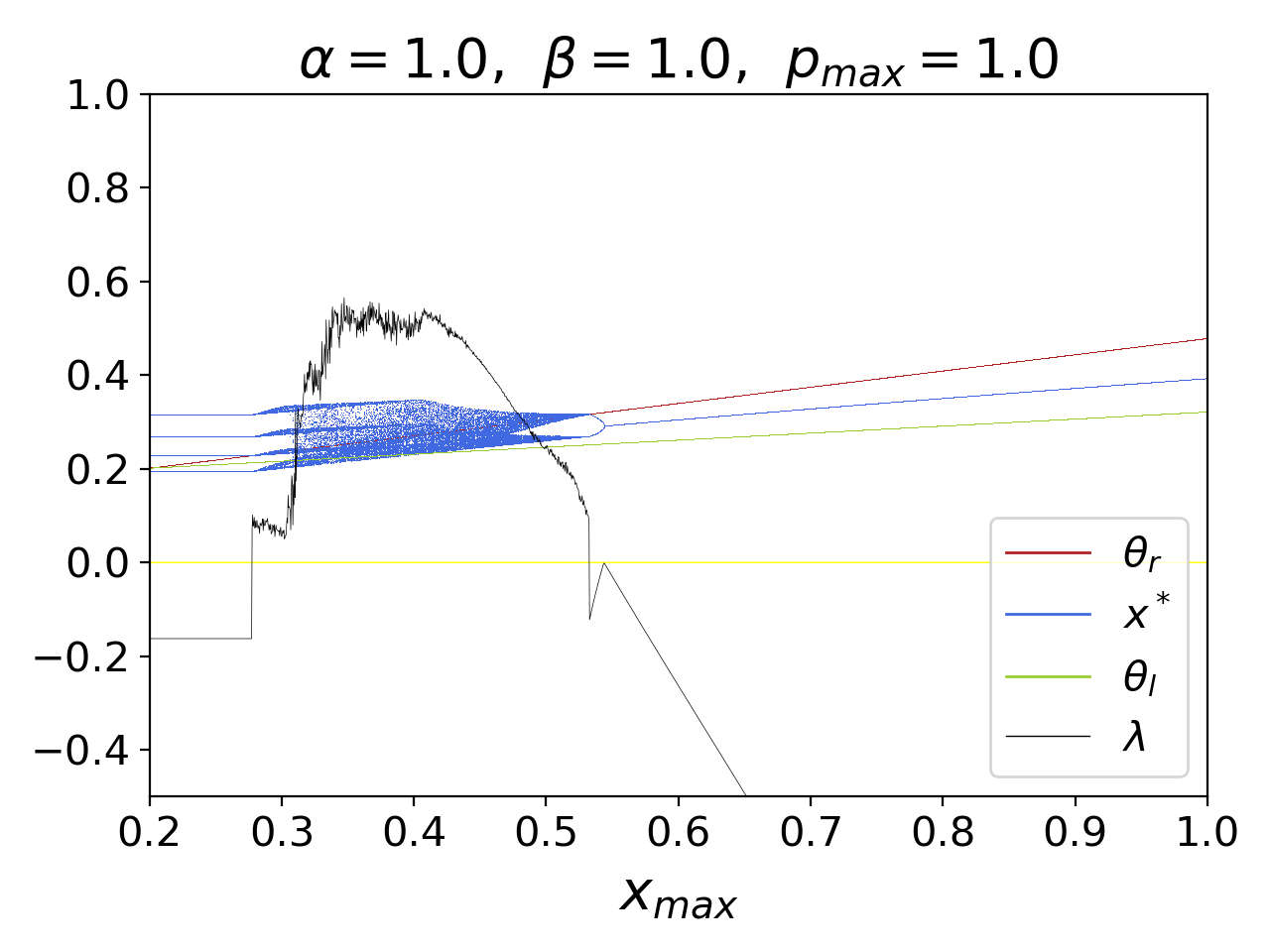} \includegraphics[width=0.45\textwidth]{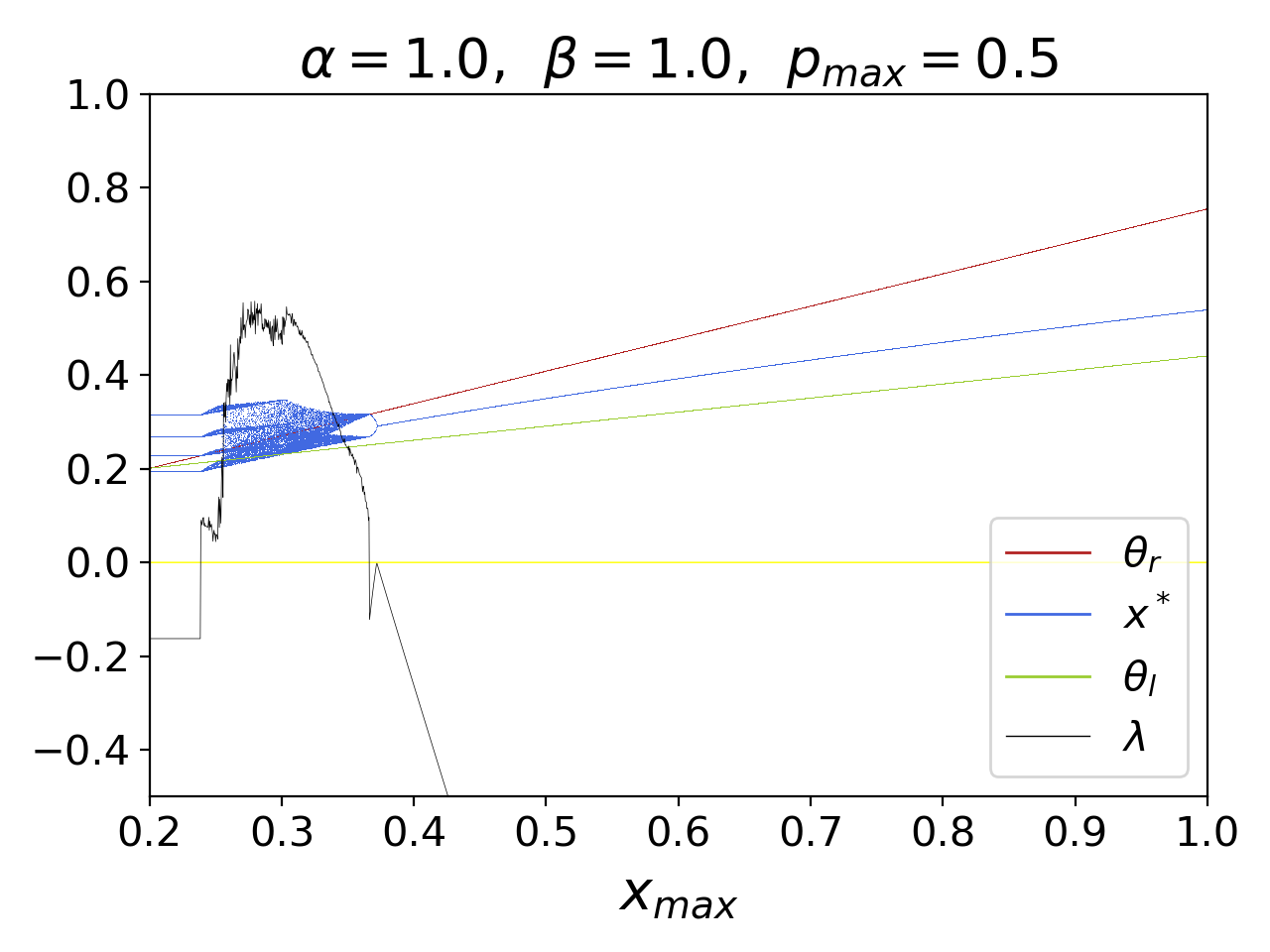} \includegraphics[width=0.45\textwidth]{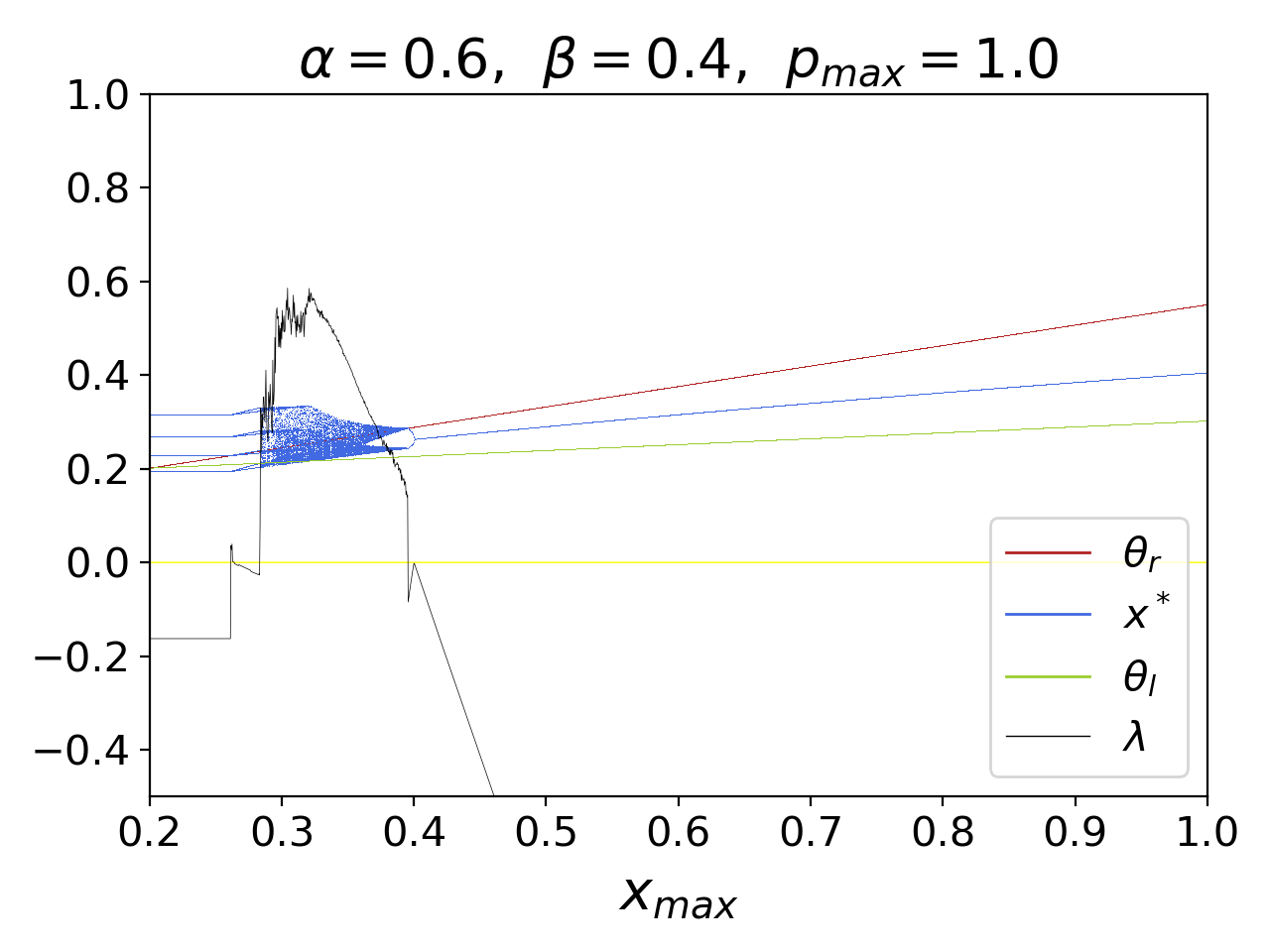} \includegraphics[width=0.45\textwidth]{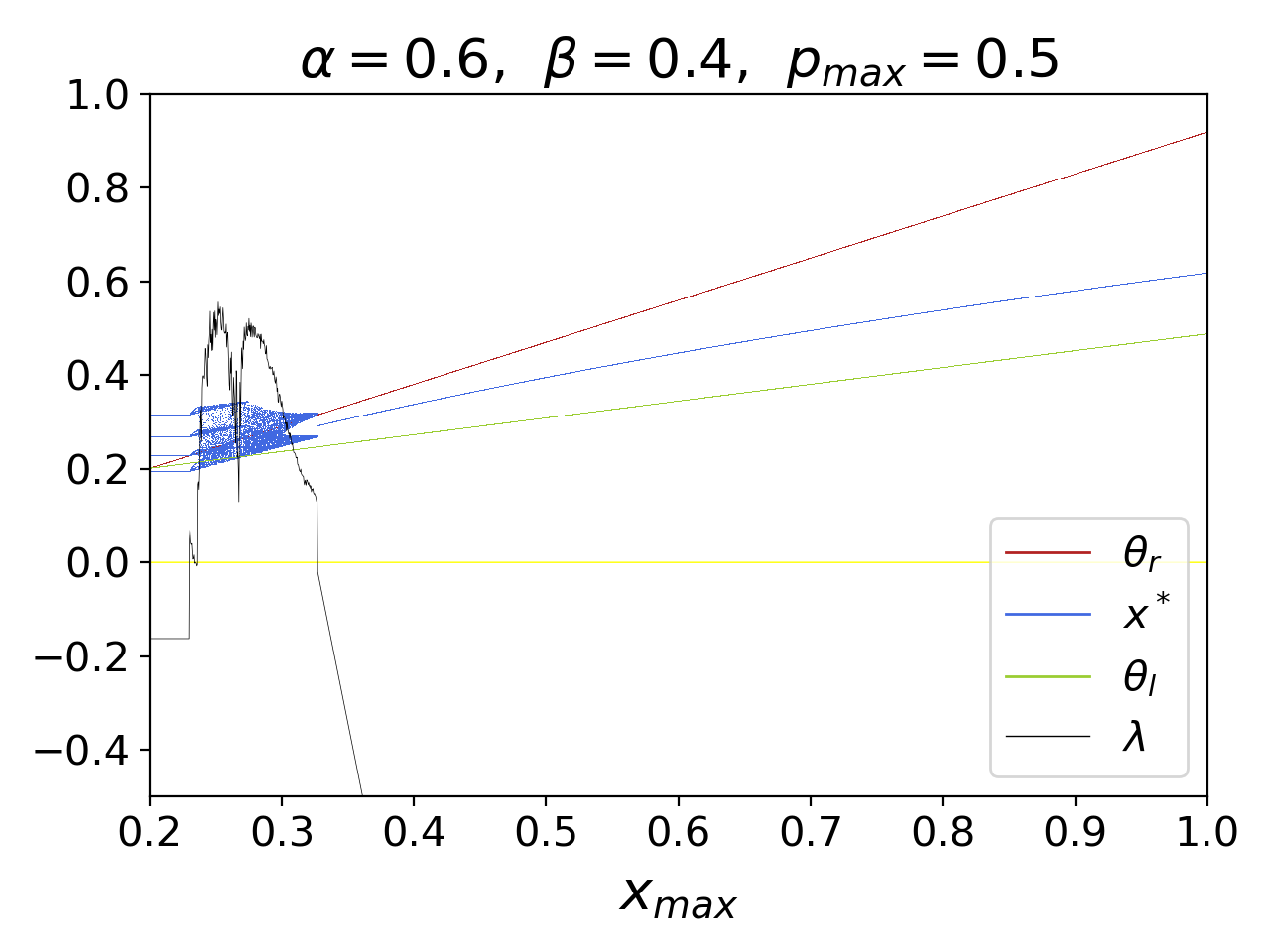}
\caption{Bifurcation diagram and Lyapunov exponents of the normalized
average queue length $x$ with respect to the parameter $x_{\max}$ for the
values of $\protect\alpha, \protect\beta$ and $p_{\max}$ displayed at the
top of the pannels. The parameter $x_{\max}$ ranges in the interval $(0.2,1]$. System parameters are given in \eqref{UMH}. Other control parameters: $x_{\min}=0.2$ and $w=0.15$.}
\label{fig-bif-xmax}
\end{figure}

Figures \ref{fig-bif-xmin} and \ref{fig-bif-xmax} show the bifurcation
diagram and the Lyapunov exponent with respect to $x_{\min }$ and $x_{\max }$, respectively, obtained with the values $\alpha =\beta =1$ (upper row) and $\alpha =0.6$, $\beta =0.4$ (lower row) for the values $p_{\max }=1$ (left
column) and $p_{\max }=0.5$ (right column). The settings of the other
control parameters are given in the captions of the figures. For clarity, we
also depict the endpoints $\theta _{l},\theta _{r}$ of the dynamical core,
that encloses the fixed point $x^{\ast }$. The parametric grid used for the
bifurcation diagrams has $2000$ points; the orbits were $550$ iterates long,
the first $500$ (the transient) having been discarded.

The most salient aspects of these figures can be summarized as follows.

(1) The bifurcations are direct for $x_{\min }$ (i.e., for growing values),
whereas they are reverse for $x_{\max }$ (i.e., for decreasing values).

(2) The bands of chaos start after the collision of an initial period-2
orbit with one of the endpoints of the dynamical core $(\theta _{l},\theta
_{r})$. One speaks in this case of a boundary collision bifurcation.

(3) The bifurcation points with respect to $x_{\min }$ occur before for $p_{\max }=1$ than for $p_{\max }=0.5$. The opposite happens with respect to $x_{\max }$. In either case we conclude that the setting $p_{\max }=0.5$ is
more robust than $p_{\max }=1.$

(4) Last but not least, the bifurcation point with respect to $x_{\min }$
(resp. $x_{\max }$) for $\alpha =0.6$, $\beta =0.4$ is greater (resp.
smaller) than for $\alpha =\beta =1$. So the model with $\alpha =0.6$, $\beta =0.4$ is more stable with respect to both $x_{\min }$ and $x_{\max }$
than the original model (\ref{Ranjan}).

Similar conclusions regarding the stability of $x^{\ast }$ were obtained in 
\cite{Amigo2020} with respect to $A_{1}$, $A_{2}$ and $w$.

\subsection{ Robustness domains}

Bifurcation diagrams are instrumental to find ranges of the control
parameters that guarantee a stable dynamics. But this is not sufficient for
our purposes. In this regard, an important feature in the design of AQM
mechanisms is that small changes in the selected control parameters should
not affect the stability of the system. In other words, the selected values
of the control parameters (especially, $\alpha $ and $\beta $) should be
also robust. As a result, robustness in the sense meant here, enables a
stable operation under both changing system parameters and perturbations due
to physical and numerical noise.

\begin{figure}[tbph]
\centering
\includegraphics[width=0.45\textwidth]{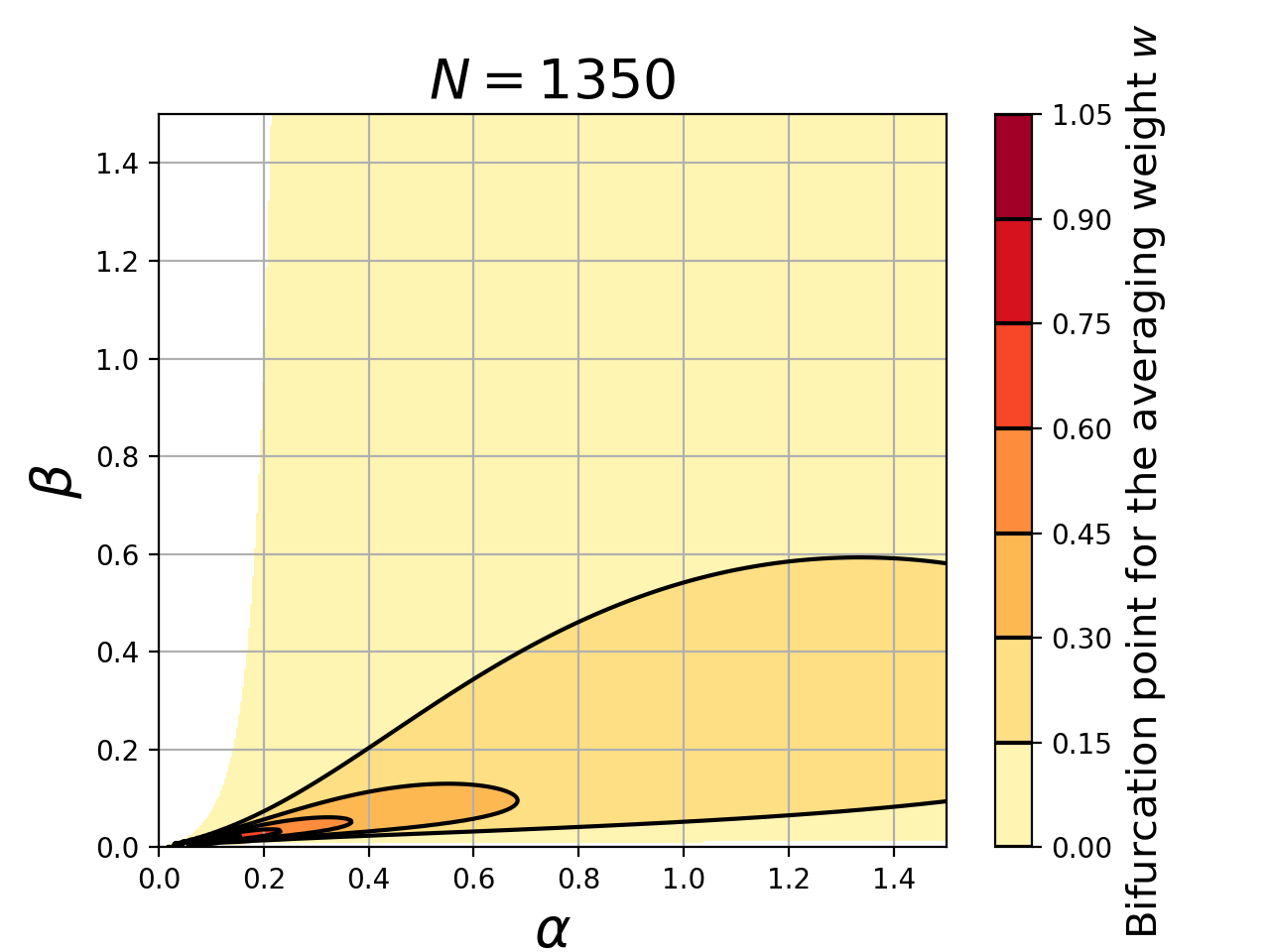} \includegraphics[width=0.45\textwidth]{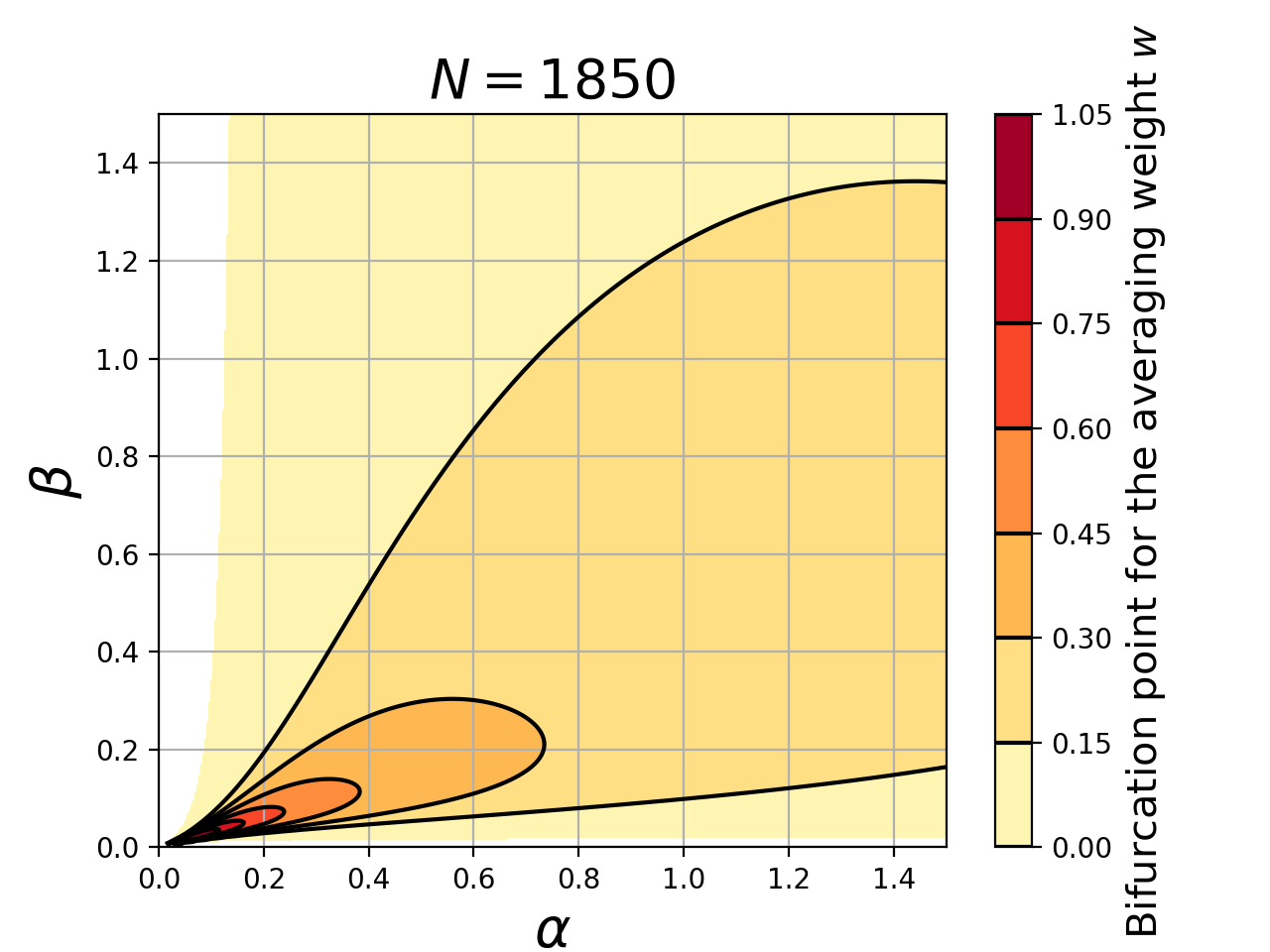}
\caption{$(\protect\alpha ,\protect\beta )$-parametric sweeps for the
bifurcation point of the averaging weight $w$. Left panel: $N=1350$; Right
panel: $N=1850$ (Source: \protect\cite{Amigo2020}). Remaining system
parameters are given in (\protect\ref{UMH}). Other control parameters: $x_{\min }=0.2$, $x_{\max }=0.6$.}
\label{fig-sweeps-wbif}
\end{figure}

If bifurcations diagrams allow to find the stability ranges of different
parameters, robustness domains allow to fine-tune the selection of control
parameters for stability control under changing system parameters. In Figure \ref{fig-sweeps-wbif} we have scanned the interval $[0.002,1.5]\times
\lbrack 0.002,1.5]$\ of the $(\alpha ,\beta )$-plane with precision $\Delta
\alpha =\Delta \beta =3.745\times 10^{-3}$\ (corresponding to a grid of $400\times 400$\ points). For each point $(\alpha ,\beta )$\ we have
calculated the bifurcation point of the averaging weight $w_{\mathrm{bif}}$\
in two different scenarios according to the number of connections: $N=1350$\
on the left panel and and $N=1850$ on the right panel.

Regarding the robustness domains shown in the left panel of Figure \ref{fig-sweeps-wbif}, we see that the second greatest domain corresponds to $0.15\leq w_{\mathrm{bif}}\leq 0.30$, and the third greatest to $0.30\leq w_{\mathrm{bif}}\leq 0.45$. Although the latter is smaller, settings of $\alpha
,\beta $ in that domain (with a typical selection $\alpha =0.5$, $\beta =0.1$) leave an ample range $0<w<0.30$ for stable operation. Points $(\alpha
,\beta )$ close to the boundaries of the robustness (same-color) domains
should be avoided. Furthermore, comparison of the left panel with the right
one shows that, as $N$ increases, corresponding domains overlap, so one can
find suitable selections of $\alpha $ and $\beta $ for $N$ ranging from $1350 $ to the reference number $1850$. The observation that the robustness
domains with $w_{\mathrm{bif}}\geq 0.15$ shrinks as $N$ decreases agrees
with the fact that a small number of connections (users) tends to disrupt
the dynamics \cite{Ranjan2004}.

Of course, a thorough analysis is necessary for the actual design of an AQM
algorithm. Our point here is that the decision on how to fix the control
parameters is a compromise between the size of the robustness domains and
stability ranges. We believe that the robustness domains are a handy tool
for that endeavour.

\section{Conclusion and outlook}

Deterministic models (mostly involving continuous time and differential
equations) are typical of the traditional sciences, in which the evolution
laws can be derived from first principles. Stochastic models are usually put
in place in the opposite case, i.e., when the laws governing a process are
unknown, hence a data-driven description is the only option left. There are
intermediate situations, though, where a deterministic description is
possible but too complex to be of any practical use. A paradigmatic example
is statistical mechanics, where the constituents move according to Newton's
laws and, yet, a probabilistic approach, where the observables are sample
averages in phase space, is the best way to bridge the gap between the
microscopic and the macroscopic scales. Other well-known example is symbolic
dynamics, which is obtained by coarse-graining the state space of a
discrete-time dynamical system.

In this paper we proceeded the other way around: the outset was RED, a
stochastic model to manage the incoming packet queue at a router's buffer.
In this specific case, the uncertainty about the arrival times of the
packets makes unfeasible a deterministic model. But if one contents oneself
with aggregate information, which behaves in a more regular way, then the
full power of a nonlinear model is available for prediction and control,
showing that stochastic and deterministic models can complement each other
advantageously. How to transit from the original RED model to a
discrete-time dynamical formulation, where the states are average queue
sizes at the buffer, was shown in Section 2. Since the resulting nonlinear
version of RED turns out to be rather sensitive to the setting of control
parameters (in particular, to the averaging weight $w$), the model was
generalized in Section 3, Equation (\ref{GDM}), by replacing the original
RED probability law (\ref{p-ave}) for packet dropping by a new one, Equation
(\ref{p_n2}), involving the beta distribution (\ref{I(x)}). This way, two
new control parameters $\alpha $ and $\beta $ are added; for $\alpha =\beta
=1$ one recovers the initial dynamical model (\ref{Ranjan}). In Section 4 we
gathered the main theoretical results on the global stability of $x^{\ast }$, the unique fixed point of the generalized RED model (\ref{GDM}) if $A_{1}<A_{2}+x_{\max }$ (Theorem \ref{Thm1}). Sufficient conditions were
given in Theorem \ref{Thm6} for a monotone $f_{(\theta _{l},\theta _{r})}$,
and in Theorems \ref{Thm7}-\ref{Thm8} for $f_{(\theta _{l},\theta _{r})}$
unimodal. Particular results for $f_{(\theta _{l},\theta _{r})}$ a $\cup $-convex function were derived in Section 5. The chaotic behavior of the RED
dynamics, proved in Section 3, was quantified by bifurcation diagrams and
Lyapunov exponents in Section 6.1, while the robustness of the fixed point
under variations of the control parameters was discussed in Section 6.2.

The ultimate goal of the results reviewed and reported in this paper is the
design of a congestion control algorithm for the data traffic on the
Internet in a realistic environment, i.e., when some parameters of the
communication network (notably, the number of users $M$ and the round-trip
time $d$) may vary over time. This was the main reason for studying convex
mappings in Section 5 and parametric robustness in Section 6.2. The former
provide an interesting trade-off between the implementation of the
analytical constraints and practical controllability; the latter is crucial
if the control parameters (especially $\alpha $, $\beta $ and $w$) have to
be adjusted in real time. Congestion control under real conditions is the
subject of work in progress.

\medskip

\section*{Acknowledgments}

This material is partially based upon work supported by the Swedish Research
Council under grant no. 2016-06596 while J.M.A. was in residence at Institut
Mittag-Leffler in Djursholm, Sweden during the summer semester 2019. This
work was financially supported by the Spanish Ministry of Economy, Industry
and Competitiveness, grant MTM2016-74921-P (AEI/FEDER, EU).


\end{document}